\documentclass [11pt]{article}
\usepackage{amsthm}
\usepackage{amssymb}
\usepackage{epsfig}

\newcounter{contador}
\newtheorem{propo}[contador]{Proposition}
\newtheorem{teo}[contador]{Theorem}
\newtheorem{lem}[contador]{Lemma}
\newtheorem{nota}[contador]{Remark}

\newtheorem{corol}[contador]{Corollary}

\newcommand{\rec}{\noindent}    
\newcommand{\dem}{\rec {\it Proof. }}  

\newcommand{\dps}{\displaystyle} 
\renewcommand{\qed}{\ \hfill\rule[-1mm]{2mm}{3.2mm}}

\newcommand{\trans}{\pitchfork}

\newcommand{\su}{{\mathbb S}^1} 

\newcommand{\esferafor}{{\mathbb S}^2} 
\newcommand{\enya}{${\rm \tilde{n}}$}
\newcommand{\bx}{\bar{x}}

\newcommand{\tD}{\tilde{\Delta}}
\newcommand{\pr}{{\rm {\bf P}}{\mathbb R}^2}
\newcommand{\ymin}{y_{\min}}
\newcommand{\lin}{{\cal L}}
\newcommand{\sg}{{\cal G}}

\newcommand{\R}{{\mathbb R}}

\newcommand{\N}{{\mathbb N}}

\title{Dynamics of the third order Lyness' difference equation
\footnote{{\bf Acknowledgements}.  The authors are grateful to
Francesc Ma\~{n}osas for his interesting and kind suggestions and
comments. GSD-UAB and CoDALab Groups are supported by the
Government of Catalonia through the SGR program. They are also
supported by DGICYT through grants MTM2005-06098-C02-01    (first
and second authors) and DPI2005-08-668-C03-1 (third author). This
research was partially done while the third author was visiting
Barcelona's CRM (Centre de Recerca Matem\`atica), having the
support the UPC's ``Programa d'ajuts a la mobilitat''; V.
Ma\~{n}osa acknowledges both institutions.}}

\author{Anna Cima$^{(1)}$, Armengol Gasull$^{(1)}$ and V\'{\i}ctor Ma\~{n}osa $^{(2)}$
  \\*[.1truecm]
{\small \textsl{$^{(1)}$ Dept. de Matem\`{a}tiques, Facultat de
Ci\`{e}ncies,}}
\\*[-.25truecm] {\small \textsl{Universitat Aut\`{o}noma de Barcelona,}}
\\*[-.25truecm] {\small \textsl{08193 Bellaterra, Barcelona, Spain}}
\\*[-.25truecm] {\small \textsl{cima@mat.uab.es, gasull@mat.uab.es}}
\\*[-.25truecm] {\small \textsl{$^{(2)}$ Dept. de Matem\`{a}tica Aplicada III,}}
\\*[-.25truecm] {\small \textsl{Control, Dynamics and Applications Group (CoDALab)}}
\\*[-.25truecm] {\small \textsl{Universitat Polit\`{e}cnica de Catalunya}}
\\*[-.25truecm] {\small \textsl{Colom 1, 08222 Terrassa, Spain}}
\\*[-.25truecm] {\small \textsl{victor.manosa@upc.edu}}}

\begin{document}

\maketitle
\begin{abstract} This paper studies the iterates of the third order Lyness'
recurrence $x_{k+3}=(a+x_{k+1}+x_{k+2})/x_k,$ with positive initial conditions,
being $a$ also a positive parameter. It is known that for $a=1$  all the
sequences generated by this recurrence are 8-periodic.  We  prove that for each
$a\ne1$ there are infinitely many initial conditions giving rise to  periodic
sequences which have  almost all the  even periods and that for a full measure
set of initial conditions the sequences generated by the recurrence are dense
in either one or two disjoint bounded intervals of $\R.$ Finally we show that
the set of initial conditions giving rise to periodic sequences of odd period
is contained in a codimension one algebraic variety (so it has zero measure)
and that for an open set of values of $a$ it also contains  all the odd
numbers, except finitely many of them.
\end{abstract}

\rec {\sl 2000 Mathematics Subject Classification:} 39A11, 39A20.

\rec {\sl Keywords:} Difference equation, discrete dynamical system, circle
map, rotation number, periodic orbit, first integral.

\break

\tableofcontents

\break

\section{Introduction and main results}

\subsection{The third order Lyness' difference equation}\label{sec11}

The excellent unpublished paper of Zeeman  \cite{Z} about the  celebrated
Lyness' second order difference equation
\begin{equation}\label{re2}
x_{n+2}=\frac{a+x_{n+1}}{x_n},\quad{\rm with}\quad a>0\,,\,x_1>0\,,\,x_2>0,
\end{equation}
gives the key points for understanding the behaviour of the
sequences generated by (\ref{re2}). In this reference it is
proved  that the map induced by (\ref{re2}),
\begin{equation}\label{fa2}
f(x,y)=\left(y,\dps{\frac{a+y}{x}}\right),
\end{equation}
defined  on $\{(x,y)\in\R^2\,:\,x>0\,,\,y>0\}$ leaves invariant
the level curves of the first integral $V(x,y)=\left( x+1 \right)
\left( y+1 \right) \left( a+x+y \right) /(xy)$ and, which is more
important, that on each set $\{V(x,y)=h\},$ the map $f$ is
conjugated to a rotation of the circle with rotation number
$\rho_a(h).$ By using this result, Zeeman explains the behavior of
all the sequences generated by  (\ref{re2}), modulus a conjecture,
the monotonous dependence of $\rho_a(h)$ with respect to $h$ once
$a\ne1$ is fixed. Recall that when $a=1$, except for the fixed
point, all the sequences generated by (\ref{re2}) are
$5-$periodic. This conjecture has been  proved to be true in
\cite{BC}. The study of the periods that can appear in the Lyness
equation, as well as the study of the rotation number has also
been done in \cite{BR}.

This paper studies a similar problem to the one considered by Zeeman but in
dimension three, and proves that in this case the dynamics are also described by
rotations. The fact that we are in a higher dimension makes the problem more
difficult.

Concretely, we consider the third order Lyness' recurrence,
\begin{equation}\label{re1}
x_{n+3}=\frac{a+x_{n+2}+x_{n+1}}{x_n},
\end{equation}
for $a>0$ and positive initial conditions $x_1,x_2$ and $x_3,$
{\it i.e.} such that $(x_1,x_2,x_3)\in
O^+:=\{(x,y,z)\in\R^3\,:\,x>0\,,\,y>0\,,\,z>0\}.$ This
recurrence  is also known as Todd's recurrence, see \cite{Grove,
KLR}. Recall that if some initial condition is such that
$(x_1,x_2,x_3)=(x_{1+p},x_{2+p},x_{3+p}),$ and $p$ is the minimal
positive number satisfying this property it is said that this
initial condition gives rise to a $p-$periodic sequence. It is
well known that when $a=1$ for any positive initial condition it
holds that all the initial conditions in $O^+$ of (\ref{re1}) are
$8,2$ or 1-periodic. We are interested to understand which is the
situation when $a\ne1.$   Our main result is:

\break

\begin{teo}\label{reclyness}
Consider the third order Lyness' recurrence (\ref{re1}) for $a>0$ and positive
initial conditions $x_1,x_2$ and $x_3$.
\begin{itemize}
\item[(i)] If $a\ne 1$ there is a computable value $q_0(a)\in\N$ such that
for any $q>q_0(a)$ there exist continua of initial conditions giving rise to
$2q$--periodic sequences.

\item[(ii)]
The set of even periods arising when $a\in(0,\infty)$ contains all the even
numbers except possibly $4,6,10,12,16,18,24,28$ and $40.$

\item[(iii)] The set of initial conditions giving rise to odd periods is contained
in an algebraic codimension one subset of $O^+.$ Moreover, there is an open set
$\mathcal U\subset(0,1)\cup(1,\infty)$  of values of $a$ for which the  set of
the odd periods contains all the odd numbers except possibly finitely many of
them.

\item[(iv)] If $a\ne 1$ then there exist a dense set  of initial
conditions in $O^+$  such that the sequence generated by  (\ref{re1}) is dense
in either one or two disjoint bounded intervals of $\R.$
\end{itemize}
\end{teo}

The above theorem makes one to wonder the following natural
questions: Are there some $a>0$ and some initial condition in
$O^+$ such that the recurrence at this initial condition is
periodic of period $4,6,10,12,16,18,24,28$ or $40$? Fixed $a>0,$
which are exactly all the even periods of the recurrence? And the
odd periods? Is the set ${\cal U}$ introduced in Theorem
\ref{reclyness} (iii) $U=(0,1)\cup(1,\infty)$?

We want to remark that when the recurrence (\ref{re1}) is
considered with initial conditions in the whole $\R^3$, the
periods that can appear can be different. For instance in
\cite{CGM1} it is proved that for some values of $a$ there are
initial (non positive) conditions giving rise to periodic
sequences with periods 2,3,4,5,6,7 and $4p$ for any $p\ge3.$

The paper is organized as follows. In  Section \ref{dds}  we
state our results on the discrete dynamical system generated by
$F$ thus obtaining the proof of Theorem \ref{reclyness}. All the
results stated in Section \ref{dds} are proved in the following
sections, moving some large proofs of the technical results  to
specific subsections and the appendices in order to improve the
readability of the paper.


\subsection{Study from a dynamical systems viewpoint.}\label{dds}

 As usual we
reduce the study of the recurrence (\ref{re1})  to the study of
the discrete dynamical system generated by the map $
F(x,y,z)=\left(y,z,(a+y+z)/{x}\right)$ defined in $O^+.$ Note that
this map is a diffeomorphism from $O^+$ to $O^+.$ A complete
description of the  discrete system generated by $F$  gives a
complete answer to the questions posed in Section \ref{sec11}, and
in particular a proof of Theorem \ref{reclyness} (see the end of
this section). Our analysis of this dynamical system is done in
two steps:

1. We will see that the phase space of $F$ is foliated by
invariant curves (sometimes degenerated to isolated points) which
are given by the level curves of two functionally independent
first integrals. The first  step is to characterize the topology
of this level sets, which turn to be diffeomorphic to circles
(when they are not isolated points).

2. The second step is  to study the dynamics of $F$ restricted to
these invariant sets.  As we will see, one of our main tools, at
this stage, will be the study of an ordinary differential equation
associated to the discrete dynamical system generated by $F$. This
is an approach different to the ones in \cite{BR}, \cite{Z} (and
even to the one in \cite{BC} although our starting point is the
same of this last reference).  Our approach turns out to be also
effective for studying other difference equations, see
\cite{CGMnou}.\newline

Nevertheless there are  some problems, named there as Questions 1
and 2, that have resisted our analysis. We remark that an answer
to them would also allow  to clarify the answers to the  questions
about~(\ref{re1}) stated in Section \ref{sec11}.

\vspace{0.5cm}

Fixed $a>0,$ consider  the diffeomorphism
\begin{equation}\label{fa}
F(x,y,z)=\left(y,z,\dps{\frac{a+y+z}{x}}\right)
\end{equation}
defined in $O^+:=\{(x,y,z)\in\R^3\,:\,x>0\,,\,y>0\,,\,z>0\}.$

We begin by introducing some  sets in $O^+$ which are invariant under the
action of $F$, in terms of the level surfaces of the well--known
(\cite{CGM2,BMG,G,K}) couple of functionally independent first integrals of
$F$, given by:
$$\begin{array}{l}
V_1(x,y,z)=\dps{\frac { \left( x+1 \right)  \left( y+1 \right) \left(
z+1 \right)  \left( a+x+y+z \right) }{xyz}},\\
V_2(x,y,z)=\dps{\frac{(1+y+z)(1+x+y)(a+x+y+z+xz)}{xyz}}.
\end{array}$$
Let $L_k=\{(x,y,z)\in O^+:\, V_1(x,y,z)=k\}$ and $M_h=\{(x,y,z)\in
O^+:\, V_2(x,y,z)=h\}$ be the level surfaces of $V_1$ and $V_2$
respectively.

The orbits of $F$ lie in  $I_{k,h}=L_k\cap M_h$ for $k\geq k_c$
and $h\geq h_c$, where $k_c$ and $h_c$ denote the values attached
at the global minima in $O^+$ of $V_1$ and $V_2$ respectively. For
a given fixed $h>h_c$, there exists $k_1=k_1(h),k_2=k_2(h)$
satisfying $k_c<k_1<k_2$ and such that $I_{k,h}\not=\emptyset$
only when $k\in[k_1,k_2].$ See Theorem \ref{conjugacio} below, or
Proposition \ref{teobifu} in Section \ref{topolevelsets} for more
details about the topology of $I_{k,h}$. We will use the notation
$A\cong B$ to denote that the two manifolds $A$ and $B$ are
diffeomorphic.

 Now, we introduce two
interesting invariant sets that will play a very important role.
The first one is $$\lin:=\{(x,(x+a)/(x-1),x)\in\R^3\mbox{ such
that }x>1\}\subset O^+.$$ It is easy to see that the set $\lin$
is a curve filled by two-periodic points of $F$ and that it
contains the unique fixed point in $O^+$: $(x_c,x_c,x_c)$, where
$x_c=1+\sqrt{1+a}.$ The second one is
$$\sg:=\{(x,y,z)\in O^+ \mbox{ such that }
G(x,y,z)=0\},$$ where
\begin{equation}\label{Ge}
G(x,y,z)=-y^3-(x+z+a+1)y^2-(x+z+a)y+xz(x+1)(z+1).
\end{equation}
The set $\lin\cup\sg$ is formed by the points in $O^+$ where the gradients of
$V_1(x,y,z)$ and $V_2(x,y,z)$ are parallel.

In particular, it is not difficult to check that \begin{equation}\label{GF}
G(F(x,y,z))=-\frac{a+y+z}{x^2}\,G(x,y,z).
\end{equation}
Note that this relation implies that $\sg$ is invariant by $F$ and
that $F$ maps the zone $\{G>0\}$ into the zone $\{G<0\}$ and
viceversa. Furthermore it implies that the dynamics of $F^2$ on
the zone $\{G>0\}$ and on $\{G<0\}$ are conjugated, being the map
$F$ itself the conjugation. Figure 1 gives an example of the more
generic position of $L_k,$ $M_h,$ $\sg$ and $\lin.$


\centerline{\includegraphics[scale=0.65]{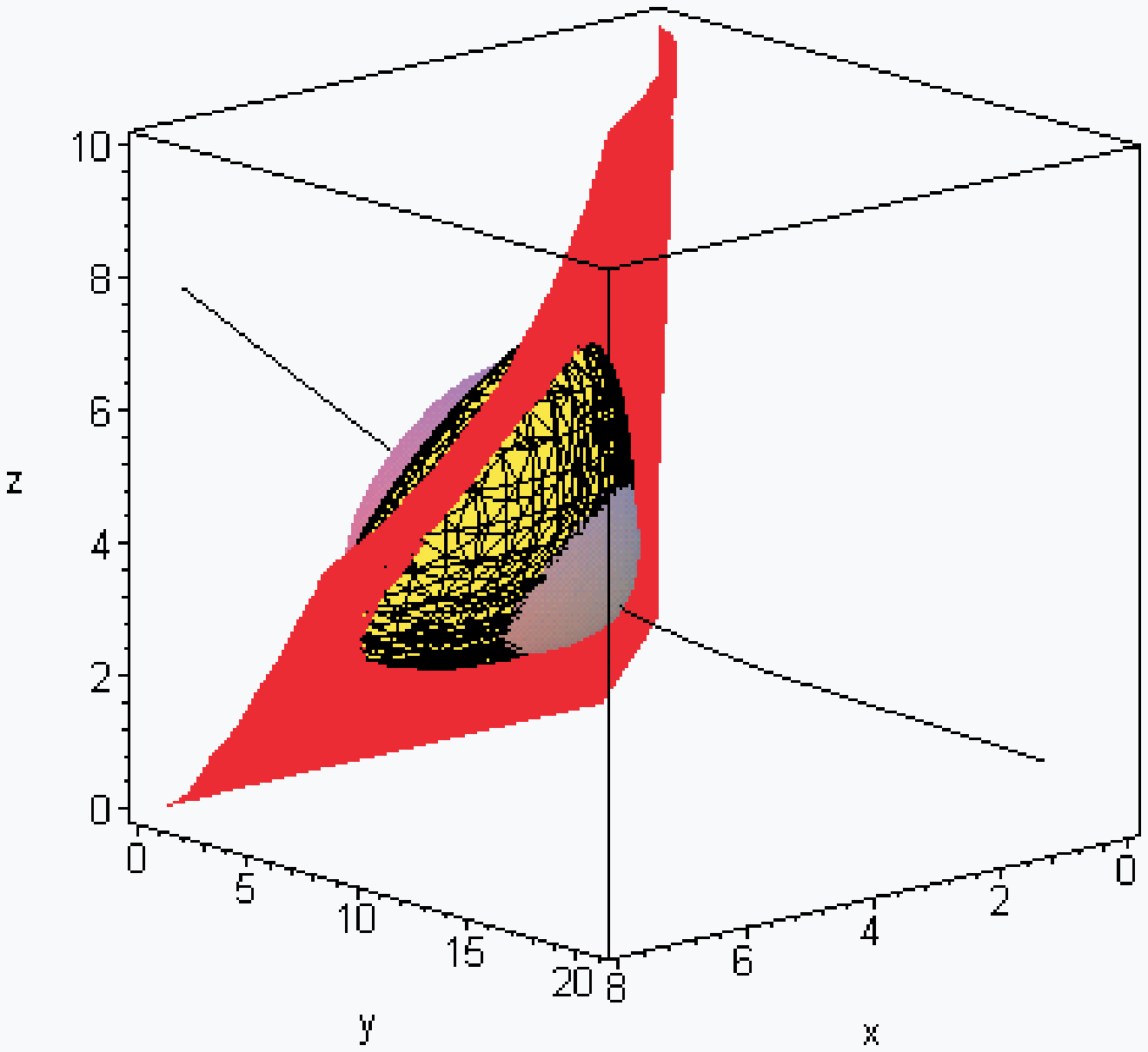}}
\centerline{Figure
1: For $a=3,$ the level surfaces $L_{31}$ and $M_{42.5}$ which are
diffeomorphic to spheres,} \centerline{ the invariant surface
$\sg$ and the line of 2-period points $\lin$.}

\vspace{0.5cm}

Our main result about the dynamics of $F$, which proved in Section
\ref{seccioconjugacio}, is:

\begin{teo}\label{conjugacio} For each fixed $h>h_c$. The following
statements hold

(i) $M_h=\cup_{k\in[k_1,k_2]} I_{k,h},$ where the  values
$k_1,k_2$ are given in Proposition \ref{teobifu}. Moreover, for
each $k\in(k_1,k_2),$ $I_{k,h}$ splits into two disjoint
connected components, that is:  $I_{k,h}=I_{k,h}^+\cup I_{k,h}^-$
where $I_{k,h}^+:=I_{k,h}\cap\{G>0\}\cong\su$ and
$I_{k,h}^-:=I_{k,h}\cap\{G<0\}\cong\su,$
$F(I^\pm_{k,h})=I^\mp_{k,h},$ $F^2(I^\pm_{k,h})=I^\pm_{k,h}$ and
the restriction of $F^2$ on each of these sets is conjugated to a
rotation of the circle with rotation number $\rho_{F^2}(k,h)$.

(ii) The set $\sg$ is invariant by $F$ and the restriction of $F$ to each set
$\sg\cap \{V_1(x,y,z)=k\}\,,\,k>k_c$ is conjugated to a rotation of the circle
with rotation number $\rho_F(k)=\frac{\rho_{F^2}(k, h(k))}2,$ where $h(k)$ is a
suitable  known function.
\end{teo}

The proof of the above result is given in two steps. The first one
is the study of the topology of the invariant sets $I_{k,h}^\pm$
(this is done in Section \ref{topolevelsets}). The second step is
to prove that over these invariant sets $F$ is conjugated to a
rotation (this is done in Section \ref{seccioconjugacio}). The
proof relies on some results that relate the rotation numbers
associated to the invariant sets $I_{k,h}^\pm$ of $F$ with the
properties of a flow constructed from $F$ which has the same
invariant sets.

Next results give some rotation numbers and periods appearing in
the dynamical system generated by $F.$ Their proofs use the
regularity of the rotation numbers varying $h,k$ and $a.$ This
regularity is studied in Section \ref{secperio}.

\begin{teo}\label{evenperiods} For  $a>0$ define
$$\rho_a:=\frac{1}{2\pi}\arccos\left(\frac{(1-a)\sqrt{1+a}}{2(1+\sqrt{1+a})
(a+1+\sqrt{1+a})}\right).$$Then for each $a\ne 1$ there are circles of initial
conditions in $\{G>0\}\setminus\lin$ and in $\{G<0\}\setminus\lin$
 such that such that $F^2$ restricted to them
is conjugated to a rotation with rotation number  taking  any value in $
\left(\frac{1}{4},\rho_a\right),$  if  $a>1,$ and any value in $
\left(\rho_a,\frac{1}{4}\right)$  if $0<a<1.$\end{teo}

In Section \ref{provamain2} we give  a constructive algorithmic approach to the
problem of determining which are all the denominators of irreducible fractions
which belong to a given interval, see Theorem \ref{racional} and Corollary
\ref{racional2}. In particular, by using these results and the above theorem,
we prove:

\begin{corol}\label{corol3} (i) For any $a\ne1$ there exists a computable value
$q_0(a)\in\N$ such that for any $q>q_0(a)$ there exist a continua of initial
conditions giving rise to $2q$--periodic orbits for $F.$

(ii) Set $$ I_{\rm rot}=\left(\frac {\pi -2\,\arcsin \left( 1/8 \right) }{4\pi
},\frac{1}{3}\right).
$$
Then, for each number $\rho$ in $I_{\rm rot}$ there exists some
$a>0$ and a circle of initial conditions such that such that $F^2$
restricted to it is conjugated to a rotation with rotation number
$\rho.$ In particular, for all the irreducible rational numbers
$p/q\in I_{\rm rot}$, there exist periodic orbits of $F^2$ of
period $q.$

  (iii) The set of even periods arising from the family
$\{F(x,y,z)=(y,z,a+y+z/x):a>0\}$ contains all the even periods except possibly
$4,6,10,12,16,18,24,28,$  and  $40.$
\end{corol}

The knowledge that we have of the odd periodic orbits of $F$ is not so detailed
as our knowledge of the even periods. We collect all our results in the
following proposition:

\begin{propo}\label{proponou} (i) All the initial conditions giving rise to
odd periods of $F$ in $O^+$ are contained in $\mathcal{G}.$

(ii) There is an open set $\mathcal{U}\subset (0,1)\cup(1,\infty)$ such that
for each $a\in \mathcal{U}$ the map $F$ over $\mathcal{G}$ has all the periods
except possibly a finite number of them.

(iii) Set $$ J_{\rm rot}=\left(\frac {\arcsin \left( 3/4 \right) }{2\pi
},\frac{1}{6}\right).
$$
For each  $\rho\in J_{\rm rot}$ there exists some $a>0$ and a circle of initial
conditions contained in $\mathcal{G}$ such that the map $F$ restricted to them
is conjugated to a rotation with rotation number  $\rho.$ Therefore, for all
the irreducible rational numbers $p/q\in J_{\rm rot}$, there exists periodic
orbits of $F$ of period  $q.$
\end{propo}

All the above results and  our numeric simulations of the functions
$\rho_{F^2}(k,k)$  and $\rho_F(k),$ detailed in  Section \ref{secnumeric}, make
as to propose the following questions. Note that the first one is similar to
Zeeman`s Conjecture.

\noindent {\bf Question 1.} \textit{ Is it true that the function
$\rho_{F^2}(k,h)$ varies monotonically when either $k$ or $h$ vary?}

 If the answer is affirmative then all the rotation numbers,  as well as the set of even periods  given in Theorem
 \ref{reclyness} and Corollary~\ref{corol3}, are the only  possible ones on $O^+\setminus\mathcal{G}.$

\noindent {\bf Question 2.} \textit{Is it true that for each
$a\neq 1$, the rotation number $\rho_{F}$ is not identically
constant on  $\sg$? Which is limit  of $\rho_{F}$ when the initial
conditions go to infinity over $\sg$?}

For instance,  the proof of Proposition \ref{proponou} (ii)
follows from the fact that, for a neighbourhood of values of
$a=\frac{3-4\cos(2\pi/7)}{(2\cos(2\pi/7)-1)^2},$  the rotation
number on $\sg$ is not constant. Unfortunately, we have not been
able to obtain a general proof of this fact. Our numerical
simulations for $a=3$ and $a=7/9$ (see Tables 1 and 3 in Section
\ref{secnumeric}) also show the same situation. If the answers to
the above questions were affirmative  we would obtain  that in
Proposition \ref{proponou}, $\mathcal {U}=(0,1)\cup(1,\infty).$

The computation of the limit of $\rho_F$ over $\cal G$ at infinity would give
us useful quantitative information about  which would be these odd periods of
$F.$

\begin{nota}\label{remark6} Note that when $a\ne1$, each map
$F$ has infinitely many different periods and has sensible
dependence with respect to the initial conditions.  This last
fact is because two close initial condition belong to two close
sets, both diffeomorphic to circles, but over each one of them
the rotation number is slightly different.
\end{nota}

\vspace{0.5cm}

{\rec {\it Proof of Theorem \ref{reclyness}.}} Parts (i) and (ii) are a direct
consequence of Corollary \ref{corol3}.  Part (iii) follows from Proposition
\ref{proponou}.

To prove (iv) observe that all initial conditions in $O^+\setminus
\lin$ give rise to rotations for $F^2$ (respectively $F$).
Moreover for most of these conditions, in the sense of Lebesgue
measure, the associated rotation numbers are irrational. Therefore
the orbits through these initial conditions are dense in a  subset
of $\cal G$ (resp. $O^+\setminus{\cal G})$ which is diffeomorphic
to $\su$ (resp. the disjoint union of two $\su$). The projection
into the $x$-axis of the orbit of $F$ coincides with the sequence
generated by (\ref{re1}). This projection is formed by one or two
disjoint closed intervals. Both situations are possible depending
if the initial conditions are near the two periodic orbit  or near
$\sg$. Hence the theorem follows. \qed

\section{Topology of the invariant sets of $F$. }\label{topolevelsets}

\subsection{The results}

This section is devoted to prove the following weaker version of
Theorem \ref{conjugacio}. Note that the difference between both
results is that in this second one the dynamics of $F$ or $F^2$ on
each of the invariant circles is not yet described. The
description of these dynamics is the goal of next sections. We use
the following notations: $A\trans B$ means that $A$ has a
transversal intersection with $B$, and $A\sqcup\,B$ means the
union of $A$ and $B$ and that both sets are disjoint. Recall also
that we say $A\cong B$ when $A$ and $B$ are two diffeomorphic
varieties.

\begin{teo}{\bf (Topology of the invariant sets)}\label{conjugaciow} Fix $h>h_c.$ Then

(i)  $M_h=\cup_{k\in[k_1,k_2]} I_{k,h},$ where the values $k_1,k_2$ are given
in Proposition \ref{teobifu}. For each $k\in(k_1,k_2),$ $I_{k,h}=I_{k,h}^+
\sqcup\, I_{k,h}^- $ and  each one of these sets is diffeomorphic to a circle.
Moreover  $F(I_{k,h}^\pm)=I_{k,h}^\mp$   and $F^2(I_{k,h}^\pm)=I_{k,h}^\pm.$

(ii)The set $\sg$ is foliated by the fix point of $F$ and  the
sets $\sg\cap \{V_1=k\}\,,\,k>k_c,$ which are invariant by $F$
and diffeomorphic to circles.
\end{teo}

The proof of the above theorem is done at the end of this
section. To prove it, we first study the level sets of $V_1$ and
$V_2$ in $O^+$ and afterwards their relative position.

Let $L_k=\{(x,y,z)\in O^+:\, V_1(x,y,z)=k\}$ and $M_h=\{(x,y,z)\in
O^+:\, V_2(x,y,z)=h\}$ be the level surfaces of $V_1$ and $V_2$
respectively. It is well known that $V_1$ has a global minimum at
$(x_c,x_c,x_c)$, where $x_c=1+\sqrt{1+a}$. We set
$$k_c=V_1(x_c,x_c,x_c)=\frac { \left( 2+\sqrt {1+a} \right) ^{3}
\left( a+3+3\,\sqrt {1+a} \right) }{ \left( 1+\sqrt {1+a} \right)
^{3}}.$$ Thus $L_k$ is not empty for $k\geq k_c$, and
$L_{k_c}=(x_c,x_c,x_c)$.

Similarly, $V_2$ also has a minimum at $(x_c,x_c,x_c)$. We set
$$h_c=V_2(x_c,x_c,x_c)={\frac { \left( 3+2\,\sqrt {1+a} \right)
^{2} \left( 2\,a+5+5\,\sqrt { 1+a} \right) }{ \left( 1+\sqrt {1+a}
\right) ^{3}}}.
$$ Then $M_h$ is not empty for $h\geq h_c$, and
$M_{h_c}=(x_c,x_c,x_c)$.

Proposition \ref{propoesferes} ( proved in Section \ref{prova8}),
states that except at the fix point all the level surfaces in
$O^+$ of $V_1$ and $V_2$ are diffeomorphic to spheres. Note that
this result proves in particular that all the orbits of $F$
starting at $O^+$ lay in  compact sets.

\begin{propo}{\bf (General properties of $L_k$ and
$M_h$)}\label{propoesferes}
\begin{itemize}
  \item[(a)] For $k>k_c,$ $L_k$ is diffeomorphic to $\esferafor.$
  \item[(b)] For $h>h_c,$ $M_h$ is diffeomorphic to $\esferafor.$
\end{itemize}
\end{propo}

Theorem \ref{conjugaciow} follows from the knowledge of the
relative positions  of the level surfaces $L_k$ and $M_h$. First
we describe the set $\mathcal F$ where $L_k$ and $M_h$ are not
transversal and give the relative position of $\mathcal F$ and
$L_k.$

\begin{lem}\label{lem1}{\bf (Locus of non-transversality of $L_k$ and $M_h$)}
Let ${\cal F}$ be the subset of $O^+$ where $\nabla V_1$ and
$\nabla V_2$ are linearly dependent, {\it i.e} ${\cal
F}:=\{\nabla V_1 \| \nabla V_2\}\cap O^+$. Then  ${\cal F}=
\lin\cup\sg$.
\end{lem}

{\rec {\it Proof.}} Some computations show that
$$\begin{array}{rl}
\left|\begin{array}{cc}
  (V_1)_x & (V_2)_x \\
  (V_1)_y & (V_2)_y
\end{array}\right|= &- \left( z+1 \right)  \left( 1+x+y \right) \left( a+z+y-xy \right) \left( a{y}^{2}+ay
-x{z}^{2}-{x}^{2}z+\right. \\
{}& \left. {y}^{2}+yz+{y}^{3}+xy-xz+{y}^{2}z-{x}^{2}{z}^{2}+x{
y}^{2} \right) /({x}^{3}{y}^{3}{z}^{2}), \end{array}$$

$$\begin{array}{rl}
\left|\begin{array}{cc}
  (V_1)_x & (V_2)_x \\
  (V_1)_z & (V_2)_z
\end{array}\right|= &- \left( y+1 \right)  \left( x-z \right) \left( a+x+y+z+xz \right) \left( a{y}^{2}+ay
-x{z}^{2}-{x}^{2}z+\right. \\
{}& \left. {y}^{2}+yz+{y}^{3}+xy-xz+{y}^{2}z-{x}^{2}{z}^{2}+x{
y}^{2} \right) /({x}^{3}{y}^{2}{z}^{3}), \qquad
\mbox{and}\end{array}$$

$$\begin{array}{rl}
\left|\begin{array}{cc}
  (V_1)_y & (V_2)_y \\
  (V_1)_z & (V_2)_z
\end{array}\right|= & \left( x+1 \right)  \left( 1+y+z \right) \left( a+x+y-yz \right) \left( a{y}^{2}+ay
-x{z}^{2}-{x}^{2}z+\right. \\
{}& \left. {y}^{2}+yz+{y}^{3}+xy-xz+{y}^{2}z-{x}^{2}{z}^{2}+x{
y}^{2} \right) /({x}^{3}{y}^{3}{z}^{2}) . \end{array}$$

The  solutions in $O^+$ of the above three functions equated to
zero satisfy either
$$ a{y}^{2}+ay -x{z}^{2}-{x}^{2}z+
{y}^{2}+yz+{y}^{3}+xy-xz+{y}^{2}z-{x}^{2}{z}^{2}+x{ y}^{2} =0,$$
which is precisely $\sg$ or $\{(x,y,z)\,:\,
  a+y-xy+z=0,\,
  x-z=0, \,
  a+x+y-yz=0\}
$ which  coincides with
$\lin=\{(x,y,z)\,:\,y=(x+a)/(x-1),\,z=x\},$ as we wanted to
prove.\qed

\vspace{0.5cm}

The topology of ${\cal F}$ is given by the next result, which is
proved in Appendix \ref{appd}.

\begin{propo}\label{lem3}{\bf (Topology of the nontransversality locus)} Fix $k>k_c$.
Then
\begin{itemize}
  \item[(i)] $\lin\cap L_k$ consists of two points which are a
  $2$--periodic orbit of $F$.
  \item[(ii)] $\sg\trans L_k$.
\item[(iii)]  $\sg\cap L_k\cong \su$.
\end{itemize}
\end{propo}

 To
describe  the relative positions  of the level surfaces $L_k$ and
$M_h$ we keep $M_h$ with $h>h_c$ fixed and consider $L_k$ for all
$k>k_c$, obtaining:

\begin{propo}{\bf (Relative positions of $L_k$ and
$M_h$)}\label{teobifu} For a given fixed $h>h_c$, there exists
$k_1:=k_1(h),k_2:=k_2(h)$ satisfying $k_c<k_1<k_2$ and such that
the following statements hold:
\begin{itemize}
  \item[(a)] If $k\in[k_c,k_1)$ then $I_{k,h}=\emptyset$.
  \item[(b)] If $k=k_1$ then  either $I_{k_1,h}=\lin \cap M_h$
  (which are two points describing a $2$--periodic orbit), or $I_{k_1,h}=\sg \cap
  M_h\cong \su$.
  \item[(c)] If $k\in(k_1,k_2)$ then $I_{k,h}\cong\su\sqcup\,\su$.
  More precisely, $I_{k,h}\cap\{G>0\}\cong\su$ and $I_{k,h}\cap\{G<0\}\cong\su.$
 \item[(d)] For  $k=k_2$ then  either  $I_{k_2,h}=\lin \cap M_h$ if  $I_{k_1,h}=\sg \cap M_h$,
 or  $I_{k_2,h}=\sg \cap   M_h\cong \su$, if $I_{k_1,h}=\lin \cap M_h$.
   \item[(e)] If $k>k_2$ then $I_{h,k}=\emptyset$.
\end{itemize}
\end{propo}

The proof of Proposition \ref{teobifu} is given in Subsection
\ref{noumeu}. Now we can prove Theorem \ref{conjugaciow}:

\vspace{0.5cm}

{\rec {\it Proof of Theorem \ref{conjugaciow}.}} The result
follows directly  from Proposition \ref{teobifu}   and the above
explained consequences  of expression (\ref{GF}).  \qed

\subsection{Proof of proposition \ref{propoesferes}}\label{prova8}

To study the surfaces $L_k,$ solving $V_1(x,y,z)=k$, we get that
they  can be written as the union of the graphs of the two
functions $z_{-}$ and $z_+$, given by:

\begin{equation}\label{zpm}
z_\pm(x,y;a,k)=\dps{\frac{\alpha(x,y;a,k)\pm\sqrt{\Delta(x,y;a,k)}}{\beta(x,y;a,k)}},
\end{equation}
defined in $\{(x,y)\in\R^2: \Delta(x,y;a,k)\geq 0\}$, where
$$ \alpha(x,y;a,k)=-a-1- \left( a+2 \right) x- \left( a+2 \right)
y-{x}^{2}-\left( a -k+3 \right) xy-{y}^{2} -{x}^{2}y -x{y}^{2},$$
$$\beta(x,y;a,k)=2(1+x+y+xy),$$ and

\rec $\Delta(x,y;a,k)=(a-1)^2+2a \left( a-1 \right) x
 + 2a\left( a-1\right)y
+ \left( 2a-2+{a}^{2} \right) {x} ^{2}\\
 + \left( -2k-2ka+ 4{a}^{2}-2
 \right) xy
 + \left( 2a-2+{a}^{2} \right) {y}^{2}
+2a{x}^{3} + \left( -4k-2ka-2+6a+2{a}^{2}
 \right) {x}^{2}y\\
+ \left( -4k-2ka -2+6a+2{a}^{2} \right) x{y}^{2} +2a{y}^{3}
+{x}^{4} + \left( 4a-2k+2 \right) {x}^{3}y\\
 + \left(
-2ka-6k+{a}^{2}+3+6a+{k}^{2}
 \right) {x}^{2}{y}^{2}
 + \left( 4a-2k+2 \right) x{y}^{3} +{y}^{4}
+2{x}^{4}y\\ + \left( -2k+4+2a \right) {x}^{3 }{y}^{2} + \left( -2k+4+2a
 \right){x}^{2} {y}^{3}
 +2x{y}^{4}
+{x}^{4}{y}^{2}+2{x}^{3}{y}^{3} +{x}^{2}{y}^{4}.$\newline

Observe that
\begin{equation}\label{observacio}
  (\alpha^2-\Delta)(x,y;a,k)=4 \left( y+1 \right) ^{2}
\left( x+1 \right) ^{2} \left( x+y+a
 \right)>0,
\end{equation}
 for $(x,y)\in Q^+:=\{(x,y):\, x>0,y>0\}$. Hence $z_\pm
(x,y;a,k)\neq 0$ on $Q^+$. This means that either $z_\pm(x,y;a,k)\in O^+$ for
all $(x,y)\in Q^+$, or $z_\pm(x,y;a,k)\in O^-:=\{(x,y,z):\, x>0,y>0,z<0\}$ for
all $(x,y)\in Q^+$. In particular, each connected component of $L_k$ with $x>0$
and $y>0$ is completely contained either in $O^+$ or in $O^-$ for all $k>k_c$.

On the other hand notice that each level surface has an equator described by
$$\displaystyle{\left.z_{\pm}(x,y;a,k)\right|_{\Delta(x,y;a,k)=0}}.$$ A description
of the planar algebraic curves $\Gamma_k:=\{(x,y)\in\R^2:\, x>0, y>0,
\Delta(x,y;a,k)=0\}$
 is given in the next lemma, which will be proved in Appendix~\ref{propoa}.
 See Figure 2 for more details. It
 is the key result to prove Proposition \ref{propoesferes} (a).

\begin{lem}\label{corbes}
For $k\geq k_c$ the planar algebraic curve
$\Gamma_k:=\{(x,y)\in\R^2:\, x>0, y>0, \Delta(x,y;a,k)=0\}$
consists of
\begin{itemize}
  \item[(a)] If $k>k_c$: two concentric ovals $\gamma_k$ and $\zeta_k$ surrounding the point
$(x_c,x_c).$ Furthermore $\gamma_k$ shrinks to $(x_c,x_c)$ when
$k\to k_c$, and if $a<1$ then the oval $\zeta_k$ has a contact
with the axis  $\{x=0\}$ and $\{y=0\}$ at $(0,1-a)$ and $(1-a,0)$
respectively.
  \item[(b)] If $k=k_c$: one oval $\zeta_{k_c}$ and the point $(x_c,x_c)\in {\rm
  Int}(\zeta_k)$.
\end{itemize}
\end{lem}


\centerline{\includegraphics[scale=0.8]{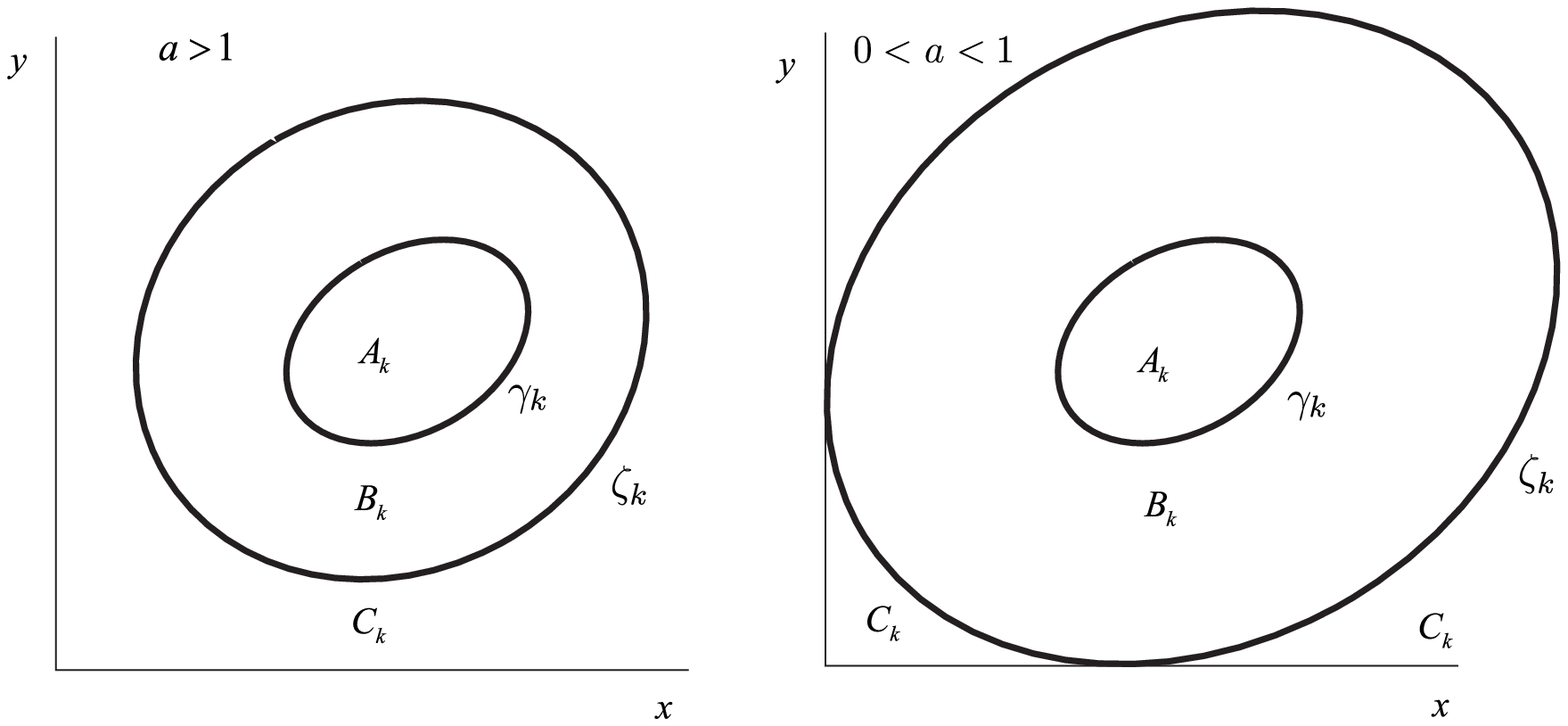}}
\centerline{Figure 2: The
 curve $\Gamma_k$ of Lemma \ref{corbes}.}

\vspace{0.5cm}

{\rec {\it Proof of Proposition \ref{propoesferes} (a). }} By
Lemma \ref{corbes}, for any $k>k_c$, $Q^+$ is split in the regions
$A_k,B_k$ and $C_k$, as is shown in Figure 2, defined in the
following way:
$$
  \begin{array}{l}
 A_k=\rm{Int}(\gamma_k),\\
 B_k=\rm{Int}(\zeta_k)\setminus\{\gamma_k\cup\rm{Int}(\gamma_k)\},\\
 C_k=Q^+\setminus\{\zeta_k\cup\rm{Int}(\zeta_k)\}.
  \end{array}
  $$

Now we will see that for any $k>k_c$  we have $\Delta(x,y;a,k)>0$
for all $(x,y)\in A_k\cup C_k$; and $\Delta(x,y;a,k)<0$ for all
$(x,y)\in B_k$. This means that the surface $\{V_1=k\}$,  defined
by (\ref{zpm}) only exists for $(x,y)\in A_k\cup C_k$.

Indeed, we can write $\Delta(x_c,x_c;a,k)=x^2y^2\, k^2+p_1(a)k+p_0(a)$, and on
the other hand $\Delta(x_c,x_c;a,k)=0$ for $k=k_1$ and $k=k_c$, such that
$k_1<k_c$. Hence, for  $k>k_c$ we have $\Delta(x_c,x_c;a,k)>0$, and this proves
that $\Delta(x,y;a,k)>0$ for all $(x,y)\in A_k$. On the other hand
$\Delta(0,0;a,k)=(1-a)^2>0$, hence $\Delta(x,y;a,k)>0$ for all $(x,y)\in C_k$.

Finally, it can be seen  that the zeros  of $\Delta(x,y;a,k)=0$ on $Q^+$ are
simple so that $\Delta(x,y;a,k)<0$ for all $(x,y)\in B_k$.

Now we observe that $\lim_{y\to 0^+}
\alpha(x,y;a,k)=\alpha(x,0;a,k)=-{x}^{2}- \left( a+2 \right)
x-a-1<0$, and that $\lim\limits_{y\to+\infty}
\alpha(x,y;a,k)=-\infty$, for all $x>0$ and $k$. This means
together with the above observation concerning equation
(\ref{observacio}), that $\{V_1=k\}\subset O^-,\mbox{ for all }
(x,y)\in C_k.$

Observe that $z_{\pm}(x_c,x_c;a,k_c)=x_c>0$, hence by continuity
$z_\pm(x_c,x_c;a,k)>0$ for $k\gtrsim k_c$. But, as seen before, by
equation (\ref{observacio}), each connected component of
$\{V_1=k\}$ with $x>0$ and $y>0$ is completely contained either in
$O^+$ or in $O^-$ for all $k>k_c$. So $z_\pm(x,y;a,k)>0$ for
$(x,y)\in A_k$.

Therefore $L_k$ is given by (\ref{zpm}) for $(x,y)\in A_k$, hence
$L_k$ is  a topological sphere. To see that  indeed it is
diffeomorphic to a sphere, by using the implicit function Theorem,
it suffices to prove that the function $V_1$ has no critical
points on $L_k.$ Computing the partial derivatives of $V_1$ we
get:
$$\begin{array}{l}
\dps{\frac{\partial V_1}{\partial x}=-\frac
{(y+1)(z+1)(-x^2+a+y+z)}{x^2yz}},\\
\dps{\frac{\partial V_1}{\partial y}=-\frac
{(x+1)(z+1)(-y^2+a+x+z)}{xy^2z}},\\
\dps{\frac{\partial V_1}{\partial z}=-\frac {(x+1)(y+1)(-z^2+a+x+y)}{xyz^2}}.
\end{array}$$
Hence the critical points of $V_1$ which lie on $O^+$ have to satisfy
$x^2=a+y+z\,,\,y^2=a+x+z\,,\,z^2=a+x+y$ which easily implies $x=y=z=x_c.$ So
the only critical point of $V_1$ on $O^+$ is the fixed point. Hence part (a) of
the proposition follows.\qed

\vspace{0.5cm}

To prove Proposition \ref{propoesferes} (b) we proceed in a
similar way that in case (a). Solving $V_2(x,y,z)=h$, we get that
the surface $\{V_2=h\}$ can we written as the union of the graph
of the two functions:

\begin{equation}\label{zpm2}
z_\pm(x,y;a,h)=\dps{\frac{\alpha(x,y;a,h)\pm\sqrt{\Delta(x,y;a,h)}}{\beta(x,y;a,h)}},
\end{equation}
where $$ \alpha(x,y;a,h)=-ya-2{x}^{2}-3x-a-{y}^{2}x-y{x}^{2}-5yx-3y+
kxy-xa-2{y}^{2}-1,$$ $$\beta(x,y;a,h)=2+2{x}^{2}+4x+2y+2yx,$$ and
$\Delta(x,y;a,h)= 1-2a+2x+2y+{h}^{2}{x}^{2}{y}^{2}-6yxa+2y{
a}^{2}x-4{y}^{2}ax-2{y}^{3}ax-4{y}^{2}a{x}^{2}-4ya{x}^{2}-2{
x}^{3}ya-6h{x}^{2}y-6kx{y}^{2}-10h{x}^{2}{y}^{2}-2h{x}^{2}{y}^
{3}-2h{x}^{3}{y}^{2}-4h{x}^{3}y-4kx{y}^{3}+4y{x}^{2}+{a}^{2}-4
xa+4yx-4ya-2kxy+{x}^{2}+{y}^{2}+4{y}^{2}x+2x{a}^{2}-2kx{
y}^{2}a-2h{x}^{2}ya-2kxya+{y}^{2}{a}^{2}+2{x}^{2}{y}^{3}+2x{y}
^{3}+2y{a}^{2}-2{x}^{2}a+2{x}^{3}y+{x}^{2}{a}^{2}-2{y}^{2}a+{y
}^{4}{x}^{2}+{y}^{2}{x}^{4}+2{x}^{3}{y}^{2}+5{x}^{2}{y}^{2}+2{y} ^{3}{x}^{3}.$

By looking at the above coefficients it is easy to check  that, if
$x>0$ and $y>0$, $\alpha(x,y;a,h)^2-\Delta>0$ and hence $z_\pm
(x,y;a,h)\neq 0$ on $Q^+:=\{(x,y):\, x>0,y>0\}$. Thus, either
$z_\pm(x,y;a,h)\in O^+$ for all $(x,y)\in Q^+$, or
$z_\pm(x,y;a,h)\in O^-:=\{(x,y,z):\, x>0,y>0,z<0\}$ for all
$(x,y)\in Q^+$, that is each connected component of $\{V_2=h\}$
with $x>0$  and $y>0$ is completely contained either in $O^+$ or
in $O^-$ for all $h>h_c$.

On the other hand observe that each level surface has an equator
given by the equation
$\displaystyle{\left.z_{\pm}(x,y;a,h)\right|_{\Delta(x,y;a,h)=0}}.$
The description of the planar algebraic curves, dropping the
subindex $a,$ $\Gamma_h:=\{(x,y)\in\R^2:\, x>0, y>0,
\Delta(x,y;a,h)=0\}$ is again the key of the proof of Proposition
\ref{propoesferes} (b), see  Figure 3. We will use the following
lemma, proved in Appendix \ref{propob}.\newline


\centerline{\includegraphics[scale=0.9]{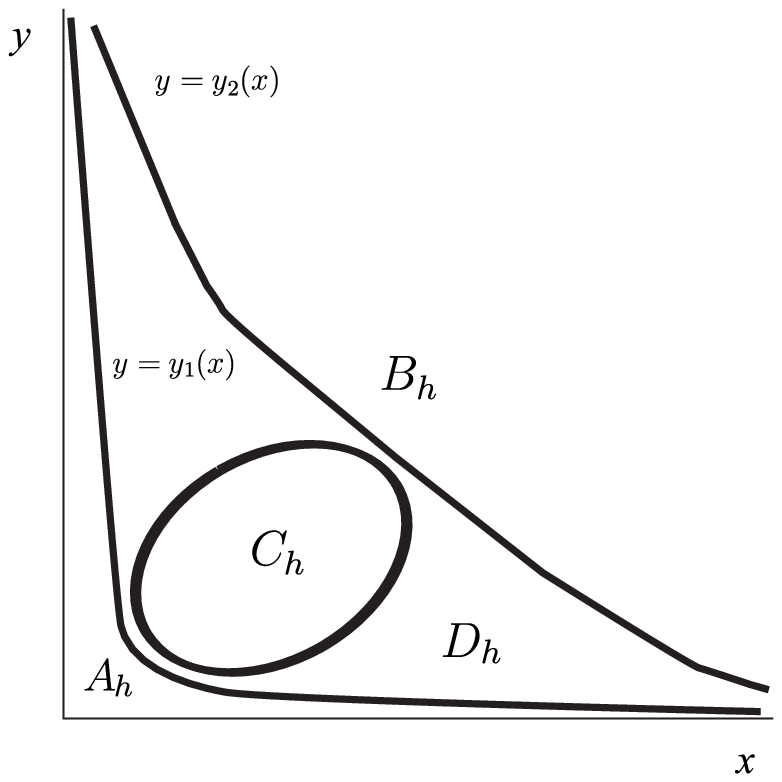}}
\centerline{Figure 3: The
 curve $\Gamma_h$ of Lemma
\ref{corbes2}.  }

\begin{lem}\label{corbes2}
For $h\geq h_c$ the planar algebraic curve
$\Gamma_h:=\{(x,y)\in\R^2:\, x>0, y>0, \Delta(x,y;a,h)=0\}$
consists of
\begin{itemize}
  \item[(a)]Two branches $y=y_i(x)$ , $i=1,2$,
 such that $y_1(x)<y_2(x)$, for $x>0$. Furthermore these two branches satisfy $\lim\limits_{x\to
0^+} y_i(x)=+\infty$ and $\lim\limits_{x\to +\infty} y_i(x)=0^+$.
  \item[(b)] An oval $\gamma_h$, contained between these to branches if
$h>h_c,$ and a single point if $h=h_c$.
\end{itemize}
\end{lem}

\vspace{0.5cm}

{\rec {\it Proof of Proposition \ref{propoesferes} (b).}} By using
Lemma \ref{corbes2} we have that for any $h>h_c,$ $Q^+$ splits in
four regions $A_h,B_h,C_h$ and $D_h$, as is shown in Figure 4,
defined in the following way:
$$
  \begin{array}{l}
 A_h=\{(x,y)\in Q^+: y\leq y_1(x)\},\\
 B_h=\{(x,y)\in Q^+: y\geq y_2(x)\},\\
 C_h=\gamma_h\cup\rm{Int}(\gamma_h),\\
 D_h=\{(x,y)\in Q^+: y_1(x)\leq y\leq
 y_2(x)\}\setminus\{\gamma_h\cup\rm{Int}(\gamma_h)\}.
  \end{array}
  $$
It is easy to check that for any $h>h_c,$ $\Delta(x,y;a,h)>0$ for
all $(x,y)\in A_h\cup B_h\cup C_h$; and $\Delta(x,y;a,h)<0$ for
all $(x,y)\in D_h$. This means that the surface $\{V_2=h\}$, is
only defined by (\ref{zpm2}) on the region $A_h\cup B_h\cup C_h$.

We now observe that $\lim_{y\to 0^+}
\alpha(x,y;a,h)=\alpha(x,0;a,h)=-1-a-(a+3)x-2x^2<0$, and
$\lim\limits_{y\to+\infty} \alpha(x,y;a,h)=-\infty$, for all $x>0$
and $h$. Hence $\{V_2=h\}\subset O^-,\mbox{ for all } (x,y)\in
A_h\cup B_h.$

We observe that $z_{\pm}(x_c,x_c;a,h_c)=x_c>0$, hence by
continuity $z_\pm(x_c,x_c;a,h)>0$ for $h\gtrsim h_c$.  But, as
seen before each connected component of  $\{V_2=h\}$ with $x>0$
and $y>0$ is completely contained either in $O^+$ or in $O^-$ for
all $h>h_c$. So $z_\pm(x,y;a,h)>0$ for  $(x,y)\in C_h$.

Therefore $M_h$ is given by (\ref{zpm2}) for $(x,y)\in C_h$, hence
$M_h$ is indeed a topological sphere. Finally, let us see that the
surface $M_h$ is a differentiable manifold for $h>h_c.$ It is
enough to see that $V_2$ has no critical points on $M_h.$ The
partial derivatives of $V_2$ are:

$$\begin{array}{l}
\frac{\partial V_2}{\partial x}=-{\frac { \left( y+z+1 \right) \left(
-{x}^{2}-{x}^{2}z+ya+{y}^{2}+yz +a+y+z \right)
}{{x}^{2}yz}},\\

\frac{\partial V_2}{\partial y}=-{\frac
{a+x+z+xza+xa+{x}^{2}+3\,xz-2\,{y}^{2}-x{y}^{2}z+za+{z}^{2}+2
\,{x}^{2}z-2\,x{y}^{2}-{y}^{2}a-2\,{y}^{2}z+2\,x{z}^{2}+{x}^{2}{z}^{2}
-2\,{y}^{3}}{x{y}^{2}z}},\\

\frac{\partial V_2}{\partial z}=-{\frac { \left( x+y+1 \right) \left(
-{z}^{2}-x{z}^{2}+ya+xy+{y}^{2} +a+x+y \right) }{xy{z}^{2}}} .
\end{array}$$

The critical points on $O^+$ have to satisfy:
$$\begin{array}{rl}
p(x,y,z):=&-{x}^{2}-{x}^{2}z+ya+{y}^{2}+yz +a+y+z  =  0\\
q(x,y,z):=& a+x+z+xza+xa+{x}^{2}+3\,xz-2\,{y}^{2}-x{y}^{2}z+za+{z}^{2}+2
\,{x}^{2}z-2\,x{y}^{2}\\
{}&-{y}^{2}a-2\,{y}^{2}z+2\,x{z}^{2}+{x}^{2}{z}^{2}
-2\,{y}^{3}  =  0\\
r(x,y,z):=&-{z}^{2}-x{z}^{2}+ya+xy+{y}^{2} +a+x+y  =  0.
\end{array}$$
Since $p(x,y)-r(x,y)=(z-x)(1+x+y+z+xz)$ we get $z=x,$ and substituting this
equality in $p(x,y,z)$ and $q(x,y,z)$ we have the system:
$$\begin{array}{l}
s(x,y)=-{x}^{2}-{x}^{3}+ya+{y}^{2}+xy+a+y+x=0,\\
t(x,y)= \left( x+y+1 \right)  \left(
{x}^{3}+3\,{x}^{2}-y{x}^{2}+xa+2\,x-2\,xy-2\,{y}^{2}+a-ya \right) =0.
\end{array}$$
Denoting  $$u(x,y)= {x}^{3}+3\,{x}^{2}-y{x}^{2}+xa+2\,x-2\,xy-2\,{y}^{2}+a-ya$$
we get
$$2s(x,y)+u(x,y)=-x^3+x^2-yx^2+xa+4x+3a+ya+2y=0$$
which let us to isolate $y$ in terms of $x:$
\begin{equation}\label{isoly}
y=\frac{x^3-x^2-(4+a)x-3a}{-x^2+a+2}.\end{equation} Then
$$s\left(x,\frac{x^3-x^2-(a+4)x-3a}{ \left( -{x}^{2}+a+2
\right) }\right)=\frac{h(x)}{\left( -{x}^{2}+a+2 \right) ^
{2}},$$ where $ h(x)=- \left( a-1+x+{x}^{2} \right)  \left(
a+2\,x-{x}^{2} \right) \left( 2\,a+xa-{x}^{3}-2-2\,x-2\,{x}^{2}
\right), $ and hence we have to consider three cases depending on
the zeros of $h(x).$

If  $a-1+x+{x}^{2}=0$ and $x>0,$ then $x=x_c$ and from (\ref{isoly}), $y=x_c.$
Since $z=x$ we get the fixed point.

If  $a+2\,x-{x}^{2}=0$ and $x>0,$ then $x=(-1+\sqrt{5-4a})/2$ and substituting
this value of $x$ at (\ref{isoly}) we see that the corresponding $y$ is
negative, so we do not need to study this case.

Now we have to consider the positive roots of
$g(x):={x}^{3}+2\,{x}^{2}+ \left( 2-a \right) x+2(1-a).$ We notice
that when $a<1$ there are not changes on the signs on the
coefficients of $g(x),$ and hence there are not positive roots of
$g(x)=0.$ When $a>1,$ then there is a unique change of signs
between the coefficients of $g(x),$ and hence,  by the Descartes
rule, and since $g(0)<0$ and $\lim\limits_{x\to +\infty}
g(x)=+\infty$, we get exactly one positive real root, say $\bar
x.$ We claim that for $x=\bar x,$ the corresponding value of $y$
given in (\ref{isoly}) is negative. To prove this observe that
$g(\sqrt{a+2})=4\,\sqrt {a+2}+4\,a+6>0,$ which implies that $\bar
x<\sqrt{a+2}.$ So, $a+2-\bar x^2>0,$ i. e., the denominator of
(\ref{isoly}) is positive. By evaluating the numerator of
(\ref{isoly}) at $\bar x,$ we get:
$$\bar x^3-\bar x^2-(4+a)\bar x-3a=\bar x^3-\bar x^2-(4+a)\bar x-3a-g(\bar x)=-(3\bar x^2+6\bar x+a+2)<0.$$
Hence, there are no critical points of $V_2$ in $O^+$ different from
$(x_c,x_c,x_c),$ and the result follows. \qed

\subsection{Proof of Proposition \ref{teobifu}}\label{noumeu}

To prove Proposition  \ref{teobifu} we need some technical
results, stated below.

\begin{lem}\label{lem2}
\begin{itemize}
  \item[(a)] $I_{k_c,h_c}=(x_c,x_c,x_c)$, and $I_{k,h_c}=\emptyset$
for $k>k_c$.
  \item[(b)] Fix $h>h_c$. Set $k_1=\min\limits_{M_h} V_1$ and
$k_2=\max\limits_{M_h} V_1$, then the following statements hold.
\begin{itemize}
  \item[(i)] $I_{k,h}\neq\emptyset$ if and only if $k\in [k_1,k_2].$
  \item[(ii)] $I_{k_i,h}\subset {\cal F}=\sg\cup\lin$, for each $i=1,2$.
\end{itemize}
\end{itemize}
\end{lem}

{\rec {\it Proof.}} (a) For $k=k_c$
$L_{k_c}=M_{h_c}=(x_c,x_c,x_c)$, hence
$I_{k_c,h_c}=(x_c,x_c,x_c)$. For $k>k_c$, since
$(x_c,x_c,x_c)\notin L_{k}$, $I_{k,h_c}=\emptyset$.

(b) Since $M_h$ is compact then there exists
$k_1=\min\limits_{M_h} V_1$ and $k_2=\max\limits_{M_h} V_1$.

Observe that $k_2\neq k_1$ because  $V_1(M_h)$ is not constant,
since otherwise \linebreak $\nabla V_1(M_h)\| \nabla V_2(M_h)$,
but this only happens in $M_h\cap(\lin\cup\sg)\neq M_h$ (to prove
this last inequality just consider that if $h>h_c$, then $M_h\cap
\lin$ consists of two points which are not contained in $\sg$,
hence $M_h\cap(\lin\cup\sg)$ is disconnected).

By definition of $k_1$ and $k_2$ it is obvious that
$I_{k,h}=\emptyset$ if $k\notin[k_1,k_2]$.

Let $a,b\in\R^3$ be a points in $M_h$ such that $V_1(a)=k_1$, and
$V_1(b)=k_2$ (observe that this points exist because the absolute
extrema of ${V_1}_{|\{M_h\}}$ are reached).

Take now a continuous curve $\gamma:[0,1]\longrightarrow M_h$
such that $\gamma(0)=a$ and $\gamma(1)=b$. The function
$$\begin{array}{rcl}
  g:[0,1] & \longrightarrow  & [k_1,k_2] \\
  t & \longrightarrow & V_1(\gamma(t))
\end{array}$$
is continuous, that is for all $k\in (k_1,k_2)$ there exists at
least  $t_k$ such that $g(t_k)=k$, and $\gamma(t_k)\in I_{k,h}$.
Hence the result follows.

(ii) By the definition of $k_1$ and $k_2$, and using the theory of
extrema with constraints, $V_1$ reaches these values  at points
where the gradient vectors of $V_1$ and $V_2$ are parallel, hence
in ${\cal F},$ as we wanted to prove.\qed

\vspace{0.5cm}

\begin{lem}\label{lem4} Fixed $k>k_c$.  $\left.{V_2}\right|_{\{L_k\cap \lin\}}=c_1$ and
$\left.{V_2}\right|_{\{L_k\cap\sg\}}=c_2$ where $c_1$ and $c_2$
are different constants.
\end{lem}

{\rec {\it Proof.}} From Lemma \ref{lem3} (i), $L_{k}\cap
\lin=p_1\sqcup\, p_2$, but since $\{p_1,p_2\}$ is a $2$--periodic
orbit and $V_2$ is an invariant
$V_2(p_1)=V_2(F(p_2))=V_2(p_2)=c_1$.

Also $V_2$ is constant over $L_{k}\cap\sg$.  Indeed, by Lemma
\ref{lem3} (ii), $L_{k}\cap \sg\cong \su$, hence
  we can consider a ${\cal C}^1$-parameterization of $L_{k}\cap \sg$
  given by $\gamma(t)$. By definition
  $\dps{\frac{d}{dt}}V_1(\gamma(t))=\nabla V_1(\gamma(t))\cdot \gamma'(t)=0$.
  But observe that on $\sg$, $\nabla V_1\| \nabla V_2$, hence
  $\nabla V_2(\gamma(t))\cdot \gamma'(t)=\nabla V_1(\gamma(t))\cdot \gamma'(t)=0$, thus
  $V_2(\gamma(t))$ is constant. Hence $V_2(L_k\cap \sg)=c_2$.

To prove that $c_1\neq c_2$, just observe that by a same argument
than the one used in the proof of Lemma \ref{lem2}, $c_1$ and
$c_2$ both must be the extrema of $V_2(L_k)$, and $V_2$ is not
constant over $L_k$.\qed

\vspace{0.5cm}

\begin{corol}\label{corol5}
\begin{itemize}
  \item[(i)] Either $$I_{k_1,h}=\lin\cap L_k\cong p_1\sqcup\, p_2,\mbox{
and }I_{k_2,h}=\sg\cap L_k\cong \su,\quad \mbox{or}$$
$$I_{k_1,h}=\sg\cap L_k\cong \su,\mbox{ and }
  I_{k_2,h}=\lin\cap L_k\cong p_1\sqcup\, p_2,$$
  where $p_1,p_2$ is the two-periodic orbit, located in $L_k$.
    \item[(ii)] If $k\in(k_1,k_2)$ then $L_k\trans M_h$. In
  particular $I_{k,h}\cong \sqcup\,_{{\rm finite}} \su$.
\end{itemize}
\end{corol}

{\rec {\it Proof.}} (i) Lemma \ref{lem2} (ii) ensures that
$I_{k_1,h}\subset {\cal F}=\lin\cup\sg$. Lemma \ref{lem4} prevents
that there exist points $a,b\in L_{k_1}\cap M_h$ such that
$a\in\lin$ and $b\in\sg$. Thus either we have

$\bullet$ Case 1: $I_{k_1,h}\subset \lin$  and hence
$I_{k_1,h}=I_{k_1,h}\cap \lin$, or

$\bullet$ Case 2: $I_{k_1,h}\subset \sg$, (hence
$I_{k_1,h}=I_{k_1,h}\cap \sg$).

Observe that by Lemma \ref{lem3} (i) $L_{k_1}\cap \lin=p_1\sqcup\,
p_2$, but since, by Lemma~\ref{lem4}, $V_2(p_1)=V_2(p_2),$  we
have that in the Case 1, $I_{k_1,h}=L_{k_1}\cap M_h\cap
\lin=L_{k_1}\cap\lin=p_1\sqcup\, p_2$.

In the second case, we need to prove that
  $I_{k_1,h}=L_{k_1}\cap M_h\cap \sg =L_{k_1}\cap \sg\cong \su$. Observe
  that this is a consequence of the fact that by Lemma \ref{lem4},
 $V_2$ is constant over $L_{k_1}\cap\sg.$

  The same argument holds if instead of $I_{k_1,h}$ we consider
  $I_{k_2,h}$. But now observe that if we are in the first case,
  then $I_{k_2,h}=L_k\cap\sg$, because each point belongs only to one level set of $V_1.$
   The same happens in the second case.

  (ii) Statement (i) implies that the locus  of non transversal intersections
  of the foliation of $O^+$ given by $\{L_k\}_{\{k>k_c\}}$ with $M_h$ are given only by
  $I_{k_i,h}$, $i=1,2$. On the other hand Lemma \ref{lem2} (i) ensures that $I_{k,h}\neq\emptyset$
  for $k\in(k_1,k_2)$, therefore $L_k$ and $M_h$ must intersect transversally in this
  case. Thus, see \cite[page 30]{GP}, ${\rm Codim}(L_k\cap M_h)={\rm Codim}(L_k)+ {\rm Codim}(M_h)=1+1=2$ and
$L_k\cap M_h$ is a submanifold of $\R^3$. This implies that
$I_{k,h}$ is a union of  curves. But since both $L_k$ and $M_h$
are compact, and each connected, compact $1$--dimensional
manifold is diffeomorphic to $\su,$ see \cite[page 208]{GP}, then
$I_{k,h}\cong \sqcup\, \su$. But these disjoint union of $\su$
lie in a compact region (say $M_h$) and are defined by analytic
equations, therefore it must be a finite union. \qed

\vspace{0.5cm}

Next lemma shows that in (ii) of the above Corollary the finite
union is exactly  two $\su.$

\begin{lem}\label{lem6}{\bf (Topology at the transversal intersections of $L_k$ and $M_h$)}
For all $k\in(k_1,k_2)$, $I_{k,h}\cong \su \sqcup\, \su$.
\end{lem}

{\rec {\it Proof.}} From Corollary \ref{corol5} (ii) and
Proposition \ref{propoesferes} (b) we know that for all
$k\in(k_1,k_2)$, $I_{k,h}\cong\sqcup\,_{{\rm finite}} \su$.

Observe  that ${\cal V}_1=\{L_k\}_{\{k\in(k_1,k_2)\}}$, induce a
foliation of closed curves  on $M_h$ nesting the $2$--periodic
points defined by $M_h\cap \lin$. These two periodic points are
the only ones in the foliation of $M_h$ induced by ${\cal V}_2=
\{L_k\}_{\{k\in[k_1,k_2]\}}$, and since they are in the plane
$z=x$, each of the closed curves of ${\cal V}_1$ must intersect
the plane $z=x$. We want to prove  there are only two of them.

Consider now the restriction of $M_h$ to the plane $z=x$, given
by the equation $V_2(x,y,x)=h$, which solutions are described by
two functions $x\rightarrow y_{\pm}(x,h)$.

Let $v_1(x):=V_1(x,y_+(x,h),x)$ be the restriction of $V_1$ over
the branch $y=y_+(x,h)$. We only need to prove that fixed
$k\in(k_1,k_2)$, the equation
\begin{equation}\label{laultima}
V_1(x,y_+(x,h),x)=k,
\end{equation}
has only two solutions, that  correspond to the two closed curves
of the statement (observe that from expression (\ref{GF}) for any
closed invariant curve $\gamma_1$ such that
$\gamma_1\cap\{G<0\}\neq \emptyset$, we have
$\gamma_1\cap\{G>0\}= \emptyset$). To see this we will prove that
the singular points of $v_1(x)$ are located in
$(\sg\cup\lin)\cap\{z=x\}$, hence $v_1$ is monotonic for those
$x$ such that $(x,y_+(x,h),x)\in \{G>0\}\cap\{z=x\}$ or
$\{G<0\}\cap\{z=x\}$, and therefore equation (\ref{laultima}) has
two solutions. Indeed, in $z=x$, $(V_i)_x=(V_i)_z$ for $i=1,2,$
\begin{equation}\label{ast}
  v_1'(x)=\left.\left((V_1)_x+(V_1)_y\dps{\frac{dy_+}{dx}}+(V_1)_z\right)\right|_{\{z=x,y=y_+\}}=
  \left.\left(2(V_1)_x+(V_1)_y\dps{\frac{dy_+}{dx}}\right)\right|_{\{z=x,y=y_+\}},
\end{equation}
and from $V_2(x,y,x)=h$, we have
\begin{equation}\label{tri}
  \frac{dy_+}{dx}=-\left.\left(\frac{(V_2)_x+(V_2)_z}{(V_2)_y}\right)\right|_{\{z=x,y=y_+\}}=
  -\left.\left(\frac{2(V_2)_x}{(V_2)_y}\right)\right|_{\{z=x,y=y_+\}}.
\end{equation}
Using equations (\ref{ast}) and (\ref{tri}) we have $v_1'(x)=0$
if and only if
$$
\frac{(V_1)_x}{(V_2)_x}=\frac{(V_1)_y}{(V_2)_y}.
$$
Hence on the locus where the gradient vectors of $V_1$ and $V_2$
are parallel. This set is $(\lin\cup\sg)\cap\{z=x\},$  as we
wanted to prove. \qed

\vspace{0.5cm}

{\rec {\it \bf Proof of Proposition \ref{teobifu}.}} Statements
(a) and (e) are a direct consequence of Lemma \ref{lem2} (i),
statements (b) and (d) from Corollary \ref{corol5} (i) and,
finally, statement (c) follows from  Lemma \ref{lem6}.\qed

\section{Dynamics of $F$. Proof of Theorem \ref{conjugacio}. }\label{seccioconjugacio}

Next result relates, under some hypotheses,  the dynamics of an ordinary
differential equation and a discrete dynamical system that share an one
dimensional invariant set.

\begin{teo}\label{teoclau}
Let $f:U\rightarrow U$ be a ${\cal C}^1$ map where $U$ is an open
connected set $U\subset \R^n$, satisfying the following
assumptions:
\begin{itemize}
\item [(A1)] There exists a $\mathcal{C}^1$ vector field $X$ in $U$ such that
\begin{equation}\label{condiciodc}
X(f(q))=(Df)_q\,X(q)\quad \it{for\,\, all}\quad q\in U.
\end{equation}
\item[(A2)] For a fixed $p\in U,$ the map $f$ leaves  invariant
$\gamma_p,$  where $\gamma_p$ is the orbit of $\dot x=X(x)$ which
passes trough $p.$  In particular there exists $\tau\in
\mathbb{R}$ such that $\varphi(\tau,p)=f(p),$ where
$\varphi(t,p)$ is the flow of $X.$
\end{itemize}
Then, if $\gamma_p\cong\su,$  the restriction of $f$ on
$\gamma_p$ is conjugated to a rotation on the circle with
rotation number $\tau/T$, where
   $T$ is the period of $\gamma_p.$
\end{teo}

\dem  By substituting $q=\varphi(t,p)$ in (\ref{condiciodc}) we obtain:
$$X(f(\varphi(t,p)))=(Df)_{\varphi(t,p)}\,X(\varphi(t,p))\quad
\mbox{for\,\, all}\quad t\in \R.$$ Notice that the above equality
precisely says that the function $t\to f(\varphi(t,p))$  is also
a solution of $\dot x=X(x).$ Since when  $t=0$ it passes trough
$f(p),$ by the theorem of uniqueness of solutions we have that
$$\varphi(t,f(p))=f(\varphi(t,p))\quad
\mbox{for\,\, all}\quad t\in I_p\cap I_{f(p)}.$$ Hence, since
$\varphi(\tau,p)=f(p),$ if $q=\varphi(t,p)$ we get that
$\varphi(t,\varphi(\tau,p))=f(q)$ or, equivalently that
\begin{equation}\label{key} \varphi(\tau,q)=f(q)\quad\mbox{ for
all}\quad q\in\gamma_p.\end{equation} Let us prove by using this
relation that the  map $f:\gamma_p\rightarrow \gamma_p$ is
conjugated to a rotation of the circle with rotation number
$\rho:=\tau/T.$

Indeed we prove that the map $h:\su\rightarrow \gamma_p$ given by
$h(\exp{(it)})=\varphi\left(\frac{T}{2\,\pi}\,t,p\right)$ is the
desired conjugation. To see this it suffices  to show that
$f\circ h=h\circ r_{\tau},$ where $r_{\tau}$ is the rotation of
angle $2\,\pi\tau/T.$ The following chains of equalities give us
the desired result.
$$(f\circ h)(\exp{(it)})\,=\,f\left(\varphi\left(\frac{T}{2\pi}t,p\right)\right)=
\varphi\left(\tau+\frac{T}{2\pi}t,p\right),$$
$$(h\circ
r_{\tau})(\exp{(it)})\,=\,h\left(\exp\left(i(t+\frac{2\pi\tau}{T})\right)\right)\,=\,
\varphi\left(\frac{T}{2\pi}(t+\frac{2\pi\tau}{T}),p\right)\,=\,
\varphi\left(\tau+\frac{T}{2\pi}t,p\right),$$ where we have used
(\ref{key}).\qed

\vspace{0.5cm}

\rec{\it Proof of Theorem \ref{conjugacio}}. By using Theorem
\ref{conjugaciow}, only  remains to  prove  that   $F$ or $F^2$ restricted to
the invariant leaves given in this theorem  are conjugated to  rotations. This
will be done by using Theorem \ref{teoclau}. To apply it we need a vector field
$X$ having the same invariant leaves that in  Theorem \ref{conjugaciow}  and
satisfying (\ref{condiciodc}). We start with
 the vector field $\tilde X=\nabla V_1\times \nabla V_2,$ where recall
that $V_1$ and $V_2$ are the invariants of $F.$ We obtain that
$$
\begin{array}{lcl}
\tilde X_1(x,y,z) :=(x+1)(1+y+z)(yz-x-y-a)G(x,y,z)/(x^2y^3z^3),\\
\tilde X_2(x,y,z) :=(y+1)(z-x)(a+x+y+z+xz)G(x,y,z)/(x^3y^2z^3),\\
\tilde X_3(x,y,z) :=(z+1)(1+x+y)(a+y+z-xy)G(x,y,z)/(x^3y^3z^2).
\end{array}
$$
Clearly it has $V_1$ and $V_2$ as first integrals, but
unfortunately it does not satisfy (\ref{condiciodc}). It is
natural to try to remove the common factors of the components of
the above vector field. We consider the new differential equation
defined by the vector field
$X(x,y,z):=\frac{(xyz)^2}{G(x,y,z)}(\nabla V_1(x,y,z)\times \nabla
V_2(x,y,z))$:

\begin{equation}\label{campcontinu}
\begin{array}{lcl}
\dot x & = &X_1(x,y,z) :=(x+1)(1+y+z)(yz-x-y-a)/(yz),\\
\dot y & = &X_2(x,y,z) :=(y+1)(z-x)(a+x+y+z+xz)/(xz)\\
\dot z & = &X_3(x,y,z) :=(z+1)(1+x+y)(a+y+z-xy)/(xy),
\end{array}
\end{equation}

A computation shows that it satisfies condition
(\ref{condiciodc}), i.e. $X(F(q))=(DF)_q\,X(q)$ in $O^+,$ and then
also $X(F^2(q))=(DF^2)_q\,X(q)$ in $O^+.$

 Since $X$ also has $V_1$ and
$V_2$ as a first integrals, each connected component of $I_{k,h}$
will be an orbit of $\dot x=X(x).$ By Theorem \ref{conjugaciow},
the sets $I_{k,h}\cap \{G>0\}\cong \su$ and $I_{k,h}\cap
\{G<0\}\cong \su,$ for $k\in(k_1,k_2),$
 are periodic orbits of $X$ and invariant by $F^2.$ Since condition
(\ref{condiciodc}) is satisfied for $F^2$, Theorem \ref{teoclau} applies and
$F^2$ is conjugated to a rotation of the circle. Hence, statement $(i)$
follows.

(ii) Now consider $\sg\cap\{V_1=k\},$ which by Theorem
\ref{conjugaciow} is also a periodic orbit of $\dot x=X(x)$ and is
invariant by $F.$  Since (\ref{condiciodc}) is also satisfied, by
using again  Theorem \ref{teoclau} we get that on
$\sg\cap\{V_1=k\},$ $F$ is conjugated to a rotation of the circle,
as we wanted to prove.\qed

\section{Properties of the rotation numbers of $F$}\label{secperio}

The  main result of this section proves  the analyticity of the rotation number
of $F^2$ in $O^+\setminus \{\lin\}$  and  computes the limit of these rotation
numbers when we tend in a certain way to the line of two periodic points
$\lin.$

\begin{propo}\label{analitiques}
(i) For each fixed $a>0$ and $h>h_c$ consider the values $k_1:=k_1(a)$ and
$k_2:=k_2(a)$ given in Theorem \ref{teobifu}. Then the assignment
$(k,h,a)\longrightarrow \rho(k,h,a),$ where $\rho(k,h,a)$ is the rotation
number of $F^2$ restricted to $I_{k,h},$ is analytic for all $(k,h,a)$ with
$a>0,$ $h>h_c$ and $k\in(k_1,k_2).$

(ii) Fix $a>0$ and  $k>k_c.$ Let $\rho_F(k)$ be the rotation number of $F$ at
each level curve $\{V_1=k\}\cap\sg,$ and let $-1, \cos(\bar\theta)\pm i
\sin(\bar\theta)$ be the three eigenvalues of $DF$ at the fix point of $F.$
Then
$$\lim_{k\to k_c}\rho_F(k)=\frac{\bar\theta}{2\pi}=\frac{1}{2\pi}
\arccos\left(\frac{a-1+\sqrt{1+a}}{2a}\right).$$

(iii)  Fix $a>0$ and consider the surface $M_h$ with $h>h_c$ a
fixed value. The set $M_h\cap\{G>0\}$ is filled by closed curves
given by $I_{k,h}^+$ and by  a fix point of $F^2$ given by
$p_h=M_h\cap\{V_1=k^*\}\cap\{G>0\}\in\lin$ given by
$p_h=(x_h,(a+x_h)/(x_h-1),x_h),\,x_h>1.$ Let $\rho_{F^2}(k)$ be
the rotation number of $F^2$ on $I_{k,h}^+$ and $1,
\cos(\theta_h)\pm i \sin(\theta_h)$ the three eigenvalues of
$DF^2$ at $p_h.$ Then,
$$\lim_{k\to k^*}\rho_{F^2}(k)=\frac{\theta_h}{2\pi}=\frac{1}{2\pi}\arccos\left(\frac{(a-1)(1-x_h)}
{2x_h(a+x_h)}\right).$$

\end{propo}

To prove the above proposition we need some preliminary results.

\begin{lem}\label{lemav1}
 (i) The only singular points of the vector field (\ref{campcontinu})  in
$O^+$ are on  ${\cal L}$. Moreover, except on ${\cal L},$ the plane
$\Sigma=\{z=x\}$ is a global transversal section for its flow.

(ii) Set $a=a_0>0,$ $h>h_c$, $k\in(k_1,k_2)$, and
$p_0=(x_0,y_0,x_0)\in\Sigma\setminus\{{\cal L}\cup{\cal G}\}$.
Then there exists a neighborhood of $p_0$ in $\Sigma$ (namely
$\Sigma_{loc}=\Sigma\cap B_{\varepsilon}(p_0)$), such that the
points $p\in \Sigma_{loc}$ depend analytically on $a, k$ and $h$.
\end{lem}

\dem A straightforward computation shows that $X$ is orthogonal to $\Sigma$
outside $\lin$ and  that ${\cal L}$ is filled by the only  singular points of
$X$ in $O^+$. So, statement $(i)$ follows.

It is important to notice that as we will see  in the proof of Lemma
\ref{lem6}, all the periodic orbits of the vector field $X$ must intersect
$\Sigma$.

To prove $(ii)$ consider
$$ \begin{array}{rlll}
  V:&\R^{5}&\longrightarrow&\R^2 \\
  {} & (x,y,a,h,k) & \longrightarrow & (V_1(x,y,x)-k,V_2(x,y,x)-h),
\end{array}
$$
where the dependence of $a$ is hidden in $V_1$ and $V_2.$ On one hand
$V(x_0,y_0,a_0,k_0,h_0)=0$, and on the other hand
$$
\det\left.\left(\begin{array}{cc}
  \frac{\partial V_1}{\partial x} & \frac{\partial V_1}{\partial y} \\
  \frac{\partial V_2}{\partial x} & \frac{\partial V_2}{\partial y}\end{array}
  \right)\right|_{z=x}=-\dps{\frac{(x+1)(x+y+1)(a+x-xy+y)\,
G(x,y,x)}{x^5y^3}},
$$
which  is nonvanishing  in $\Sigma\setminus\{{\cal L}\cup{\cal G}\}$. From the
implicit function theorem, $x$ and $y$ are analytic functions of $a,k$ and $h$
in a neighborhood of $p_0=(x_0,y_0,x_0)$, namely $\Sigma_{loc}$.\qed

\vspace{0.5cm}

\begin{lem}\label{lemav3} Let $\varphi(t,p,a)$ be the flow of (\ref{campcontinu}),
where we  explicitly write the dependence with respect to $a.$
 Fix $a>0$ and take $p =(x,y,x)\in\Sigma\setminus\{{\cal
L}\cup{\cal G}\}$. \begin{itemize}
\item [(i)] If $T(p,a)$ is the period of the periodic orbit
of $X$ passing through $p$, then $T(p,a)$ is analytic at $(p,a).$
\item [(ii)] If $\tau(p,a)$ is defined by the equation
$\varphi(\tau(p,a),p,a)=F^2(p),$  then $\tau(p,a)$ is also analytic at $(p,a)$.
\end{itemize}
\end{lem}

\dem Consider the system $\varphi(t,p,a)-p={0}.$
 Obviously
$\varphi(T(p,a),p,a)-p={0}$, and  $$ \dps{\frac{\partial
}{\partial
t}\left(\varphi(t,p,a)-p\right)}=X(p)=(X_1(p),0,-X_1(p))\neq {0},
$$ because $X_1(p)\neq 0$ in $\Sigma\setminus\lin.$ By applying the implicit
function theorem to the first component  of the above system we have that in a
neighborhood of $(p,a)$ in $\Sigma\setminus\{{\cal L}\cup{\cal G}\}$ the period
function $T(p,a)$ is analytic. The proof of $(ii)$ follows applying the same
argument to equation $\varphi(t,p,a)-F^2(p)={0}$. \qed

\vspace{0.5cm}

\begin{lem}\label{difeocondiciodc} Suppose that we have a smooth vector field
$\dot{q}=X(q)$ and a smooth map $f$ in a neighborhood
${\cal{U}}\subseteq\R^n$, satisfying condition (\ref{condiciodc}).
If $p=h(q)$, and $h$ is a diffeomorphism between ${\cal{U}}$ and
$h(\cal{U})$, then the induced vector field
$\tilde{X}(p):=(Dh)_qX(q)$ and the map $\tilde{f}=h\circ f\circ
h^{-1}$ also satisfies condition (\ref{condiciodc}).
\end{lem}

\dem Indeed
$\tilde{X}(p):=(Dh)_qX(q)=(Dh)_{h^{-1}(p)}X\left(h^{-1}(p)\right)$.
Hence
$$
\begin{array}{ll}
\tilde{X}\left(\tilde{f}(p)\right) & =(Dh)_{h^{-1}(h(f(h^{-1}(p))))}X\left(h^{-1}(h(f(h^{-1}(p))))\right)= \\
 {} &=(Dh)_{f(h^{-1}(q))}\left(Df\right)_{h^{-1}(p)}X\left(h^{-1}(p)\right)=\\
 {}
 &=\left[(Dh)_{f(h^{-1}(p))}\left(Df\right)_{h^{-1}(p)}(Dh^{-1})_p\right]\cdot\left[(Dh)_{h^{-1}(p)}X\left(h^{-1}(p)\right)\right]=\\
{} &=D\tilde{f}(p)\cdot\tilde{X}(p).\end{array}
$$\qed

\begin{lem}\label{formanormal} Consider the planar vector field
\begin{equation}\label{FN}
X(u,v):=-v\,g(u^2+v^2)\dps{\frac{\partial}{\partial
u}}+u\,g(u^2+v^2)\dps{\frac{\partial}{\partial v}},\end{equation}
with $g(0)\ne0$ and $f$ a differentiable map in a neighborhood of
the origin, such that it leaves invariant the circles
$\gamma_r:=\{u^2+v^2=r^2\}$ and  condition (\ref{condiciodc}) is
satisfied, {\it i.e.} $X(f(u,v))=(Df)_{(u,v)}X(u,v).$ Then $f$ is
conjugated on each $\gamma_r$ to a  rotation with rotation number
$\rho(r)$ and $\lim\limits_{r\rightarrow 0}\rho(r)=\theta, $
where $\cos \theta\pm i\sin\theta$ are the eigenvalues of
$\left.Df\right|_{(0,0)}$.
\end{lem}

\dem
 By  Theorem \ref{teoclau}, and taking
 into account that the sets $\gamma_r$ are also invariant under $X$,
 we obtain that on each $\gamma_r,$ $f$ is conjugated to a rotation and there exists $\tau(\sqrt{u^2+v^2})$, such that
 $$f(u,v)=\varphi\left(\tau(\sqrt{u^2+v^2}),(u,v)\right)=\varphi\left(\tau(r),(u,v)\right),$$ where
 $\varphi$ is the flow of $X$. By taking polar coordinates it is not difficult
 to obtain that
$$
 \varphi(t,(u,v))=\left(\begin{array}{cc}
  \cos\left(g(r^2)t\right) & -\sin\left(g(r^2)t\right) \\
 \sin\left(g(r^2)t\right) & \cos\left(g(r^2)t\right)
\end{array}\right)\left(\begin{array}{c}
  u \\
  v
\end{array}\right) .
$$

 Hence on each set $\gamma_r$,
 $f$ is indeed the rotation of angle $\alpha(r)=g(r^2)\tau(r),$ which has
 rotation number $\rho(r)=\alpha(r)/(2\pi).$  By using this fact
 and  the differentiability of $f$
 at the origin we have
 $$
f(u,v)=\left(\begin{array}{cc}
  \cos \alpha(r) & -\sin \alpha(r) \\
 \sin \alpha(r) & \cos \alpha(r)
\end{array}\right)\left(\begin{array}{c}
  u \\
  v
\end{array}\right)=\left(\begin{array}{cc}
  \cos \alpha(0) & -\sin \alpha(0) \\
 \sin \alpha(0) & \cos \alpha(0)
\end{array}\right)\left(\begin{array}{c}
  u \\
  v
\end{array}\right)+O_2(u,v).
 $$
Hence $\alpha(0)=\theta$ and therefore $\lim_{r\to 0}\rho(r)=\lim_{r\to
0}\alpha(r)/(2\pi)=\theta/(2\pi),$ as we wanted to prove. \qed

\vspace{0.5cm}

\rec {\it Proof of Proposition \ref{analitiques}.} (i) From the
above lemmas we know that the functions
$\tau(k,h,a):=\tau(p(k,h),a)$, and $T(k,h,a):=T(p(k,h),a)$ are
analytic functions in $\Sigma_{loc}$. Since by Theorem
\ref{teoclau} we know that $\rho(k,h,a)=\tau(k,h,a)/T(k,h,a),$
then the rotation number is analytic as well.

(ii) Consider the map $F$ restricted on $\sg$ Since
$$G(x,y,z)=x(x+1)z^2+(x(x+1)-y(y+1))z-y^3-(1+a+x)y^2-(a+x)y,$$
then equation $G=0$ is equivalent to $z=z_{\pm}(x,y)$ where
$$z_{\pm}(x,y)=\frac{y(y+1)-x(x+1)\pm\sqrt{\Delta(x,y)}}{2x(x+1)}$$
and $\Delta(x,y)= (y(y+1)-x(x+1))^2+4x(x+1)(y^3+(1+a+x)y^2+(a+x)y).$ If $x>0$
and $y>0$ then $z_+(x,y)>0$ and $z_-(x,y)<0.$ Consequently the surface $\sg$
can be described as:
$$\sg=\{(x,y,z_+(x,y)):x>0,y>0\}.$$
Hence, in a neighborhood of the fix point, $\left.F\right|_{\sg}$
can be though as the planar map $ \bar F(x,y)=(y,z_+(x,y))$ in
$U=\{(x,y)\in\R^2:x>0,y>0\}.$

Clearly the map $\bar F(x,y)$ has $(x_c,x_c)$ as a fix point and the matrix
$\left.D\bar F\right|_{(x_c,y_c)}$ has the eigenvalues given by
$\lambda=\cos\bar\theta\pm i\sin\bar\theta$ where
$$\bar\theta=\arccos\left(\frac{1+x_c}{2x_c}\right)=
\arccos\left(\frac{a-1+\sqrt{1+a}}{2a}\right).$$

Let $X(x,y,z)$ be the vector field  given by (\ref{campcontinu}).
Then the map $\bar F(x,y)$ has an associated vector field
$$
\left.X\right|_{{\sg}}=:\bar
X(x,y)=X_1(x,y,z_{+}(x,y))\dps{\frac{\partial}{\partial x}}
+X_2(x,y,z_{+}(x,y))\dps{\frac{\partial}{\partial y}},
$$
 which is the restriction of the vector field
(\ref{campcontinu}) on $\sg.$ It can be checked that the vector
field $\bar X$ has the point $(x_c,x_c)$ as a singular point,
which is a non-degenerated center and satisfies
$\bar{X}(\bar{F})=(D\bar{F}) \bar{X}.$ Via an analytic change of
variables, $\bar{X}$ is conjugated with a vector field given in
the normal form (\ref{FN}), say $\bar X_N$. Through this
conjugation we also obtain that $\bar{F}$ is conjugated with a new
map $\bar{F}_N$ (which by Lemma~\ref{difeocondiciodc} satisfies
condition (\ref{condiciodc}) with $\bar X_N$). The maps $\bar F$
and $\bar{F}_N$ share the same eigenvalues at their respective
fixed points. Using Lemma \ref{formanormal} we have that
$\lim_{k\to k_c}\rho_F(k)=\bar\theta/(2\pi),$ as we wanted to
prove.

(iii)    The study of this case  is similar to one of (ii), where here $F$ and
the surface $\sg$ are replaced by $F^2$ and  $M_h,$ respectively.

A tedious computation shows that the characteristic polynomial of $(DX)_{p_h}$
is given by $P(\lambda)=\lambda\left(\lambda^2+p(x_h,a)\right)$, where
$$
p(x_h,a)={\frac { \left( x_h+1 \right)  \left( 2x_h+a-1 \right) \left(
-a+3x_ha- x_h+1+2{x_h}^{2} \right)  \left( {x_h}^{2}+x_h+a-1 \right) ^{2}}{
\left( a+x_h
 \right) ^{2}{x_h}^{2} \left( x_h-1 \right) ^{2}}}.$$
Since $x_h>1$ and $a> 0$, we have that $p(x_h;a)>0$, hence the
eigenvalues of $(DX)_{p_h}$ are 0 (which corresponds to the
tangential direction of $\lin$), and a couple of conjugated pure
imaginary ones. By the implicit function theorem, in a
neighbourhood of $p_h$ the set $M_h^+:=M_h\cap\{G>0\}$ is a
differentiable manifold of dimension 2, invariant by $X.$ Thus $X$
restricted to $M_h^+$ induces a two dimensional vector field
having a non-degenerated center at $p_h$. At this point the proof
follows in the same manner than in (ii). The computation of
$\theta_h$ is straightforward.\qed

\vspace{0.5cm}

The following result will be  useful to study the odd periods of $F.$ Although
it seems natural that it is true for any $a\ne1,$ we have not been able to
provide a general proof.

\begin{propo}\label{senar}
Consider $a=a^*:=\frac{3-4c}{(2c-1)^2}\simeq 8.29590,$ where
$c=\cos(2\pi/7).$ Then, there exists $\varepsilon>0$ such that
for any value of $a$ satisfying $|a-a^*|<\varepsilon,$ the
rotation number of $F$ over the invariant curves $\{V_1=k\}$
which foliate $\cal{G}$ is not constant.\end{propo}

\dem First we prove the result for $a=a^*.$ We proceed by
contradiction. If the set of rotation numbers were degenerated to
a point, by Proposition \ref{analitiques} this value should be the
value of the limit when we tend to the fix point, which is
$$
\frac{\arccos \left( {\frac {a^*-1+\sqrt {1+a^*}}{2a^*}} \right)}{2\pi}=\frac 1
7.
$$
Indeed we have chosen $a^*$ to obtain this value. It gives  the smallest
denominator of all the rational numbers given by the expression $(\arccos
\left( {\frac {a-1+\sqrt {1+a}}{2a}} \right))/(2\pi),$ where $a\in(0,\infty).$
By Theorem \ref{conjugacio} we also would have that $F^7$ restricted to
$\cal{G}$ would be the identity.

On the other hand, take  another point in $\cal{G},$ for instance
$q=(1,1,\frac{\sqrt{8c^2-12c+5}}{2c-1})\simeq (1,1,3.20872).$ To
prove that $F^7(q)\ne q $ it is convenient for the moment to
consider $F$ with $a=\frac{3-4d}{(2d-1)^2}$ and
$r=(1,1,\frac{\sqrt{8d^2-12d+5}}{2d-1}),$ being $d$  an unknown
parameter. The equation that forces that the first components of
$F^7(r)$ and $r$ coincide is
$${\textstyle\frac{ \left( 51-272 d +540 {d}^{2}-464 {d}^{3}+144 {d}^{4} \right)\sqrt{8d^2-12d+5} -
110+724 d-1900 {d}^{2}+2504 {d}^{3}-1664 {d}^{4}+448 {d}^{5}}{
\left(38 -168 d +292 {d}^{2} -240 {d}^{3}+80 {d}^{4} \right)
\sqrt{8d^2-12d+5}-75 +470 d- 1156 {d}^{2}+1400 {d}^{3}-832
{d}^{4}+192 {d}^{5}}=1. }$$  Working with the above equation we
obtain that its solutions are included in the solutions of
$$
64(d-1)^2(d^2-d/2-1/8)(d^2-3d/2+5/8)(128d^4-64d^3-128d^2+104d-19)=0.
 $$
Since the value $d=c$ is not a solution,  we have got that
$F^7(q)\ne q,$ which is in contradiction with our initial
assumption. Thus for $a=a^*$ we have proved that the set of
rotation numbers on $\cal G$ is not degenerated to a point. Recall
that in Proposition \ref{analitiques} (i) it is proved that the
rotation number varies continuously with respect to initial
conditions and the parameter $a.$ From this result we obtain that
the  set of all the rotation numbers over $\cal G$  is not
degenerated to a point for the values of $a$ in some neighbourhood
of $a^*,$ as we wanted to prove.  Notice that when $a=1$ this
rotation number over $\cal G,$ and also over $O^+\setminus\lin,$
is reduced to the value $1/8.$ \qed

\section{On the set of periods of $F$. Proof of Theorem \ref{evenperiods}}\label{provamain2}

In this section we prove Theorem \ref{evenperiods} and its
consequences (Corollary \ref{corol3} and Proposition
\ref{proponou}). Firstly,  we present a constructive way for
obtaining the denominators of the irreducible fractions which
belong to a given interval.

\begin{teo}\label{racional} Fix a real open interval $I=(a,b)$ and denote by $p_1=2,p_2=3, p_3,\ldots,
p_n,\ldots$ the set of all the prime numbers, ordered following
the usual order. Associated to $I$ we consider the following
natural numbers:
\begin{itemize}
\item[(i)] The smallest prime number  $p_{m+1}$ satisfying that $p_{m+1}>\max(3/(b-a),2),$
\item[(ii)] Given any prime number $p_n,$ $1\le n\le m,$ the smallest natural number $s_n$
such that $p_n^{s_n}>4/(b-a).$
\end{itemize}
By using the above numbers, define the following finite subset
of~$\mathbb{N}:$
$$F_{s_1,s_2,\ldots, s_m}:=\{n\in \mathbb{N} \,:\, n=p_1^{t_1} p_2^{t_2}\cdots
p_m^{t_m}\,\mbox{ with }\quad 0\le t_i\le s_i-1,\,i=1,2,\ldots,
m\}.$$ Then for any $r\in \mathbb{N}\setminus F_{s_1,s_2,\ldots,
s_m}$ there exists and irreducible fraction $q/r$ such that
$q/r\in I.$
\end{teo}

Next result  easily follows from the above Theorem:

\begin{corol}\label{racional2} Fix an open real interval $(a,b).$ Following the notations of the above theorem
consider the number $p:=p_1^{s_1-1} p_2^{s_2-1}\cdots
p_m^{s_m-1}.$ Then, for any $r>p$ there exists an irreducible
fraction $q/r$ such that $q/r \in (a,b).$
\end{corol}

\vspace{0.5cm}

{\rec {\it Proof of Theorem  \ref{racional}.}} We prove the
following two assertions:
\begin{itemize}
\item [(a)] If $p$ is a prime number and $p\ge p_{m+1}$ then for any natural
number $k\ge1$ there exists an irreducible fraction of the form
$\displaystyle\frac{q}{kp}\in I.$
\item[(b)] If $p_i$ is any prime number  $p_i< p_{m+1}$ and $s_i$ is the integer number given in the
statement of the theorem,  then for any natural number $k\ge1$
there exists an irreducible fraction of the form $\displaystyle
\frac{q}{kp_i^{s_i}}\in I.$
\end{itemize}

 Clearly the theorem follows from them.

Let us prove the first one. From the fact that $p_{m+1}>3/(b-a),$
we have that if $p$ is a prime number and $p\ge p_{m+1}$ then
there exists an $\ell$ such that
\begin{equation}\label{ap-1}
a<\frac{\ell-1}p<\frac\ell p <\frac{\ell+1} p <b,
\end{equation}
where the three fractions are irreducible. Hence we have proved
our assertion (a) for $k=1.$ Take now any $k>1.$ From the above
inequalities we have that
$$a<\frac{k\ell-k}{kp}<\frac{k\ell-1}{kp}<\frac{k\ell} {kp} <\frac{k\ell+1} {kp} <\frac{k\ell+k}{kp}<b.$$
Note that either $\displaystyle \frac{k\ell-1}{kp}$ or
$\displaystyle \frac{k\ell+1}{kp}$   have to be irreducible
because the factors of $k$ never divides their numerators and if
both were reducible the number $p$ should divide both numbers
$k\ell\pm1.$ Taking their difference we would have that $p$
divides $2,$ a contradiction. Thus assertion  (a) is proved.

Let us prove assertion (b). Fix any prime number $p=p_n,$ smaller
that $p_{m+1}$ and consider its associated number $s=s_n.$ From
the inequality $p_n^{s_n}>4/(b-a)$ we have that

$$a<\frac{j-1}{p^s}<\frac j {p^s} <\frac{j+1} {p^s} <\frac{j+2}{p^s}<b,$$

Note that either $j+1$ or $j$ have to be coprime with $p$ hence
taking $\ell$ either $j$ or $j+1$ we have that
$$a<\frac{\ell-1}{p^s}<\frac\ell {p^s} <\frac{\ell+1} {p^s} <b,$$
being the fraction $\ell/ p^s$ irreducible, like in (\ref{ap-1}).
When $p>2$ we can argue as in the previous case and assert that
either $\displaystyle \frac{k\ell-1}{kp^s}$ or $\displaystyle
\frac{k\ell+1}{kp^s}$ have to be irreducible, proving our result.
When $p=2$ we consider
$$a<\frac j {2^s} <\frac{j+1} {2^s} <\frac{j+2}{2^s}<b.$$
Taking $k>2$ we have that
$$a<\frac {kj} {k2^s} <\frac{kj+1} {k2^s} <\frac{kj+2}{k2^s}< \frac{kj+k}{k2^s}< b,$$
and again one of the fractions $\displaystyle \frac{kj+1} {k2^s},$
$\displaystyle \frac{kj+2} {k2^s},$ has to be irreducible, as we
wanted to prove.

\noindent{ }\qed

\vspace{0.5cm}

In the sequel we  prove  Theorem \ref{evenperiods}, Corollary \ref{corol3} and
Proposition \ref{proponou}.

\vspace{0.5cm}

\rec {\it Proof of Theorem \ref{evenperiods}}. For each $a>0,$ $a\ne 1$ and
each $x>1$ consider the function
$$r(x)=\frac{1}{2\pi}\arccos\left(\frac{(a-1)(1-x)}{2(ax+x^2)}\right).$$
Recall that from Proposition \ref{analitiques} (iii),  the
function $r(x)$ gives the limit of $\rho_{F^2}(k)$ when $k$ tend
to $V_1(p_x)$, where $p_x$ is the point on ${\cal L}$ given by
$(x,\frac{a+x}{x-1},x) $. Observe that $r(x)$ has a unique
critical point which is a maximum (resp. a minimum) at
$x=x_c=1+\sqrt{1+a}$ when $a>1$ (resp. $0<a<1$). Furthermore
$r(1)=1/4$ and $\lim_{x\to\infty}r(x)= 1/4.$ Now consider the
value $r(x_c)$ and denote it by $\rho_a:$
\begin{equation}\label{roa}
\rho_a=\frac{1}{2\pi}\arccos\left(\frac{(1-a)\sqrt{1+a}}{2(1+\sqrt{1+a})(1+a+\sqrt{1+a})}
 \right).
\end{equation}
Take $a>1$ and a number $\rho^*\in(1/4,\rho_a)$ (the case $a<1$
and $\rho^*\in(\rho_a,1/4)$ can be studied in a similar way). Let
us see that there is a continuum of initial conditions in $\{
G>0\}\setminus\cal{L}$ such that their rotation number is $\rho^*$
(and notice that by using expression (\ref{GF}), the images by $F$
of these initial conditions satisfy the same property and are in
$\{ G <0\}\setminus\cal{L}$). For $\varepsilon>0$ small enough,
there are two periodic points of $F$ in $\cal{L},$  say
$p^{\pm}=(x^\pm,\frac{a+x^\pm}{x^\pm-1},x^\pm),$  such that
$r(p^\pm)= \rho^*\pm\varepsilon.$ By Proposition \ref{analitiques}
(i) there exist initial conditions
$r^\pm\in\{V_1=V_1(p^{\pm})\}\cap\{G>0\}$ such that their
respective rotation numbers,  $\varrho^{\pm}$ satisfy
$\rho^*-2\varepsilon<\varrho^-<\rho^*<\varrho^+<\rho^*+2\varepsilon.$
Joining $r^-$ and $r^+$ by a continuous path $\Gamma\subset\{
G>0\}\setminus\cal{L},$ and by using again the continuous
dependence of the rotation number  with respect the initial
conditions, we obtain the existence of a point $r\in\Gamma$ such
that its rotation number is exactly $\rho^*.$  By
Theorem~\ref{conjugacio} the same happens with all the points  in
$\{G>0\}$ of $\{V_1=V_1(r)\}\cap\{V_2=V_2(r)\}\cong\su,$ as we
wanted to prove.\qed

\vspace{0.5cm}

\rec {\it  Proof of Corollary \ref{corol3}}. (i) By using Theorem
\ref{evenperiods} and Corollary \ref{racional2} the result follows.

(ii) Observe that the function $\rho_a$ given in (\ref{roa}) is an increasing
function such that
$$\lim_{a\to 0^+} \rho_a={\frac {\pi -2\,\arcsin \left( 1/8 \right) }{4\pi }}\simeq
0.23005\quad \mbox{ and }\quad\lim_{a\to +\infty} \rho_a=\frac{1}{3}.$$
Therefore, by using again Theorem \ref{evenperiods}, for each number in $((\pi
-2\,\arcsin \left( 1/8 \right))/(4\pi),1/4)$ there exists some $a\in(0,1)$ and
some initial condition outside $\cal G$ with this rotation number for $F^2$.
Similarly, for each number in $[1/4,1/3)$ there exist some $a\geq 1$ and some
initial condition, also outside $\cal G,$  with this rotation number. In
particular, for all the irreducible rational numbers $p/q$ with the property
$$\frac {\pi -2\,\arcsin \left( 1/8 \right) }{4\pi }<\frac{p}{q}<\frac{1}{3}$$
we can find  a value of $a$ such that $F^2$ has continua of periodic orbits  of
period $q.$

(iii) Setting $a=(\pi -2\,\arcsin \left( 1/8 \right) )/(4\pi)$, $b=1/3$ and
using the notation introduced in Theorem \ref{racional}, we have that, $m=10$,
$p_{11}=31$ and $p_1=2$ (with $s_1=5$), $p_2=3$ (with $s_2=4$), $p_3=5$ (with
$s_3=3$), $p_4=7$,  $p_5=11$, $p_6=13$, $p_7=17$, $p_8=19$, $p_9=23$ and
$p_{10}=29$ (where $s_i=2$ for $i\in\{4,\ldots,10\}$). From Theorem
\ref{racional}, we have that for all $q\in \N$, such that
  $q>q_0:=2^4\cdot3^3\cdot5^2\cdot 7\cdot 11\cdot 13\cdot 17\cdot
  19\cdot 23\cdot 29=2\,329\,089\,562\,800$ there exists some $a>0$ and
  some $q$--periodic orbit for $F^2_a$. It is now easy to
  develop a finite  algorithm in order to find which irreducible
  fractions $p/q$ with $q\le q_0$ are in $I_{\rm rot}.$  Implementing this algorithm
we get that there appear irreducible fractions with all  the denominators
except $1,2,3,5,6,8,9,12,14$ and $20.$ Doubling these numbers, and taking into
account that $\cal L$ is full of two periodic points of $F,$  we obtain (iii).
\qed

\vspace{0.5cm}

\rec {\it  Proof of Proposition \ref{proponou}}.  (i) This result is a direct
consequence of expression (\ref{GF}).

(ii) Fix a value of $a$ of the ones given in Proposition
\ref{senar}.  By this value  we know that the set of all rotation
numbers of all the points of $\cal G$ contains an open interval.
By applying Corollary \ref{racional2}  to this interval the result
follows.

(iii) Similarly that  in the proof  of (i) of Theorem
\ref{evenperiods}, for each $a>0,$ we consider the function
$s(a)=\frac{1}{2\pi}\arccos\left(\frac{a-1+\sqrt{1+a}}{2a}\right).$
Recall that from Proposition \ref{analitiques} (iii),  this
function  gives the limit of the rotation numbers over $\cal G$
when  we approach to the fix point. The range of this function
when $a>0$ is $J_{\rm rot}.$ Taking into account the continuity of
the rotation number with respect to initial conditions and the
parameter $a$, and arguing as in the last part of proof of Theorem
\ref{evenperiods},  the result follows. \qed

\section{Some numerical results}\label{secnumeric}

In this section we present some numerical explorations which lead
us to establish  the open questions stated in Section \ref{dds}.

The following tables of rotation numbers have been obtained using
the relation (\ref{key}), by numerical integration of the vector
field (\ref{campcontinu}) using a 7-8$^{\mbox{th}}$ order
Runge--Kutta method. Table~1 has been obtained taking $a=3$, and
gives the rotation number associated to the orbit passing through
some points of the surface ${\cal G}$. These points have been
taken by considering the following path over  ${\cal G}$:
\begin{equation}\label{primercami}
p(t)=\left(x(t),y(t),z(t)\right)=\left(x_c+t,\frac{a+x(t)}{x(t)-1},z(x(t),y(t))\right),
\end{equation}
where $z(x(t),y(t))$ is one of the two branches of solutions of
equation $G(x(t),y(t),z)=0$. Recall that $\tau(p):=\tau_F(p)$ is
given by the relation $\varphi(\tau_F(p),p)=F(p),$  $T(p)$ is the
period of the periodic orbit of (\ref{campcontinu}) passing
through $p$ and the rotation number of $F$ at the orbit starting
at $p$ is $\rho_F(p)=\tau_F(p)/T(p).$ Similarly we can define
$\tau_{F^2}(p)$  and $\rho_{F^2}(p).$ Note also that when both
numbers have sense $\rho_{F^2}(p)=2\rho_F(p).$ In general we
observe that the function $t\rightarrow \rho(x(t),y(t),z(t))$
seems to be decreasing.

\begin{center}
{\small
\begin{tabular}{|c|c|c|c|c|c|c|}
\hline
  $t$ & Point $p$ & $T(p)$ & $\tau_F(p) $ & $\rho_F(p) $ & $\rho_{F^2}(p) $ \\
\hline \hline
  $0  $ & $(3,3,3)$ & $--$ & $-- $ & $0.13386 $ & $0.26772 $\\
\hline
  $1 $ & $\left(4, 7/3, 1.62395  \right) $ & $0.41781 $ & $0.05586  $ & $0.13369  $ & $0.26737  $ \\
\hline
  $2 $ & $\left(5, 2, 1.06969 \right) $ & $0.36063  $ & $0.04810  $ & $0.13337  $ & $0.26674  $ \\
\hline
  $3 $ & $\left(6, 9/5, 0.78049  \right) $ & $0.30622  $ & $0.04074  $ & $0.13305  $ & $0.26610 $ \\
\hline
  $4 $ & $\left(7, 5/3, 0.60637  \right) $ & $0.26009  $ & $0.03452  $ & $0.13274 $ & $0.26549 $ \\
\hline
  $5 $ & $\left(8, 11/7, 0.49153  \right) $ & $0.22226  $ & $0.02944  $ & $0.13247  $ & $0.26494  $ \\
\hline
  $6 $ & $\left(9, 3/2, 0.41083  \right) $ & $0.19148  $ & $0.02532  $ & $0.13223 $ & $0.26446 $ \\
\hline
  $7 $ & $\left(10, 13/9, 0.35143 \right) $ & $0.16635  $ & $0.02196 $ & $0.13201  $ & $0.26402  $ \\
\hline
  $8 $ & $\left(11, 7/5, 0.30610  \right) $ & $0.14570  $ & $0.01921  $ & $0.13182  $ & $0.26364  $ \\
\hline
  $9 $ & $\left(12, 15/11, 0.27051  \right) $ & $0.12860  $ & $0.01693  $ & $0.13164  $ & $0.26328  $ \\
\hline
  $10 $ & $\left(13, 4/3, 0.24191 \right) $ & $0.11430  $ & $0.01503  $ & $0.13149  $ & $0.26298  $ \\
\hline
  $11 $ & $\left(14, 17/13, 0.21849 \right) $ & $0.10225  $ & $0.01343  $ & $0.13134  $ & $0.26269  $ \\
\hline
  $12 $ & $\left(15, 9/7, 0.19899  \right) $ & $0.09202  $ & $0.01207  $ & $0.13122 $ & $0.26244 $ \\
\hline
  $13 $ & $\left(16, 19/15, 0.18253 \right) $ & $0.08325  $ & $0.01091  $ & $0.13110  $ & $0.26220  $ \\
\hline
\end{tabular}
}
\end{center}

\vspace{0.5cm}

\centerline{Table 1. Rotation number on $\mathcal G$ for $a=3.$}

\vspace{0.5cm}

Note that the  results of Table 1  also give light to know which
are the odd periods for $F$ for a given value of $a.$ For
instance when $a=3$ it seems clear that $0.1333\ldots=2/15$ is
one of the rotation numbers reached by $F$ over $\cal G.$ Hence
for this value of $a,$ $F$ must have periodic points of period
15. By applying a three dimensional Newton method to the system
$F^{15}(x,y,z)=(x,y,z)$ we have obtained the approximated solution
$r\simeq(2.00557,5.20647, 9.89389).$  Note that $|F^{15}(r)-r|_1<
3.8 \times 10^{-5},$ where as usual $|(x,y,z)|_1=|x|+|y|+|z|.$
Indeed  there should exist infinitely many 15-periodic points,
given by all the orbits starting at the periodic orbit of
(\ref{campcontinu}) passing through  a given 15-periodic point.
We also have checked that
$$\left| F^{15\times 10^n}(r)-r \right |_1< 2.8\times 10^{n-6}\quad\mbox{for}\quad n=1,2,3,4.$$

If $r$ were a  true 15-periodic point the above values should
have been zero, but as we have already noticed in Remark
\ref{remark6} the dynamical system generated by $F$ has sensible
dependence with respect to initial conditions.

Table 2 is again obtained taking $a=3$. Now the rotation number of $F^2$ is
computed for some points in the curve given by
\begin{equation}\label{segoncami}
p(t)=\left(x(t),y(t),z(t)\right)=\left((x_c+9)\cdot(1-t)+t\,x_1,\frac{a+x(t)}{x(t)-1},z(x(t),y(t))\right),
\quad t\in[0,1],
\end{equation}
where  $z(t)$ is a fixed branch of the two branches of solutions of equation
(\ref{zpm}) and  $x(1)=x_1$ is the first coordinate of the two periodic point
of $F,$
 $p(1)\simeq\left(1.11929, 34.53097, 1.11929\right)$. This
curve joints the point
 $p(0)=\left(12, 15/11, 0.27050\ldots\right)\in{\cal G}$ with  $p(1)$
running over the level surface $$V_1=k^*:={\frac {28561}{43560}}\,{\frac {
\left( 19+3\,\sqrt {89} \right)
 \left( 197+\sqrt {89} \right) }{-79+13\,\sqrt {89}}}
 \simeq 146.70452.$$

\vspace{0.5cm}

\begin{center}

\noindent{\small
\begin{tabular}{|c|c|c|c|c|}
\hline
 $t$ & Point $p$ & $T(p)$ & $\tau_{F^2}(p) $  & $\rho_{F^2}(p) $  \\
\hline \hline
  $0  $ & $\left(12, 1.36364, 0.27051 \right) $ &  $0.12860 $ &  $0.03386 $ & $0.26328 $   \\
\hline
  $0.1  $ & $\left(10.91193, 1.40355, 0.24737 \right) $ &  $0.12857 $ &  $0.03385 $ & $0.26327 $   \\
\hline
  $0.2  $ & $\left(9.82386, 1.45332, 0.22528 \right) $ &  $0.12848 $ &  $0.03382 $ & $0.26325 $   \\
\hline
  $0.3  $ & $\left(8.73579, 1.51708, 0.20426 \right) $ &  $0.12829 $ &  $0.03376 $ & $0.26320 $  \\
\hline
  $0.4  $ & $\left(7.64772, 1.60171, 0.18437 \right) $ &  $0.12796 $ &  $0.03367 $ & $0.26310 $   \\
\hline
  $0.5  $ & $\left(6.55965, 1.71947, 0.16576 \right) $ &  $0.12743 $ &  $0.03351 $ & $0.26295 $   \\
\hline
  $0.6  $ & $\left(5.47158, 1.89454, 0.14879 \right) $ &  $0.12655 $ &  $0.03324 $ & $0.26270 $   \\
\hline
  $0.7  $ & $\left(4.38350, 2.18221, 0.13432 \right) $ &  $0.12501 $ &  $0.03279 $ & $0.26226 $  \\
\hline
  $0.8  $ & $\left(3.29543, 2.74259, 0.12498 \right) $ &  $0.12205 $ &  $0.03191 $ & $0.26145 $  \\
\hline
  $0.9  $ & $\left(2.20736, 4.31300, 0.13351 \right) $ &  $0.11498$ &  $0.02985 $ & $0.25962 $  \\
\hline
  $0.95  $ & $\left(1.66333, 7.03020, 0.17364 \right) $ &  $0.10648 $ &  $0.02744 $ & $0.25768$  \\
\hline
  $0.99 $ & $\left(1.22810, 18.53618, 0.43340\right) $ &  $0.09080 $ &  $0.02314 $ & $0.25484 $  \\
\hline
  $0.999  $ & $\left(1.13017, 31.72824, 0.95253 \right) $ &  $0.08589 $ &  $0.02183 $ & $0.25414 $  \\
\hline
  $0.9999  $ & $\left(1.12038, 34.22789, 1.09981 \right) $ &  $0.08575 $ &  $0.02179 $ & $0.25412$  \\
\hline
  $1  $ & $\left(1.11929, 34.53097, 1.11929 \right) $ &  $-- $ &  $-- $ & $0.25412 $  \\
\hline
\end{tabular}}

\end{center}
 \centerline{Table 2. Rotation number on $\{V_1=k^*\simeq 146.70452\}$ when $a=3.$}

\vspace{0.5cm}

Table 3 has been obtained by repeating the first experiment but taking $a=7/9$.
So it  gives the rotation number associated to the orbit passing through some
points of the surface ${\cal G}$, by considering the path given by
(\ref{primercami}). Notice that in this case the rotation number seems an
increasing function of $t.$

\begin{center}

{\small
\begin{tabular}{|c|c|c|c|c|c|c|}
\hline
  $t$ & Point $p$ & $T(p)$ & $\tau_F(p) $ & $\rho_F(p) $ & $\rho_{F^2}(p) $  \\
\hline \hline
$0  $ & $(7/3,7/3,7/3)$ & $--$ & $--$ & $0.12338 $ & $0.24676$ \\
\hline
$1$ & $(10/3, 37/21, 1.11361 )$ & $0.48969 $ & $0.06043 $ & $0.12340 $ & $0.24681 $ \\
\hline
$2$ & $(13/3, 23/15, 0.71973 )$ & $0.39978 $ & $0.04935 $ & $0.12345 $ & $0.24690 $ \\
\hline
$3$ & $(16/3, 55/39, 0.52965 )$ & $0.32588 $ & $0.04025 $ & $0.12350 $ & $0.24700$ \\
\hline
$4$ & $(19/3, 4/3, 0.41853 )$ & $0.26908 $ & $0.03324 $ & $0.12354 $ & $0.24708 $ \\
\hline
$5$ & $(22/3, 73/57, 0.34583)$ & $0.22552 $ & $0.02787 $ & $0.12358 $ & $0.24716 $ \\
\hline
$6$ & $(25/3, 41/33, 0.29462 )$ & $0.19171 $ & $0.02370 $ & $0.12361 $ & $0.24722 $ \\
\hline
$7$ & $(28/3, 91/75, 0.25662 )$ & $0.16502 $ & $0.02040 $ & $0.12364 $ & $0.24729 $ \\
\hline
$8$ & $(31/3, 25/21, 0.22731 )$ & $0.14363$ & $0.01776 $ & $0.12367$ & $0.24734 $ \\
\hline
$9$ & $(34/3, 109/93, 0.20401 )$ & $0.12622 $ & $0.01561 $ & $0.12369 $ & $0.24739 $ \\
\hline
$10$ & $(37/3, 59/51, 0.18506 )$ & $0.11187 $ & $0.01384 $ & $0.12372 $ & $0.24744 $ \\
\hline
$11$ & $(40/3, 127/111, 0.16933 )$ & $0.09990 $ & $0.01236 $ & $0.12374 $ & $0.24748$ \\
\hline
$12$ & $(43/3, 17/15, 0.15607 )$ & $0.08980 $ & $0.01111 $ & $0.12376 $ & $0.24752 $ \\
\hline
\end{tabular} }
\end{center}
\centerline{Table 3. Rotation number on $\mathcal G$ for $a=7/9.$}

Finally, Table 4 is obtained taking $a=7/9$, considering the path
(\ref{segoncami}), which runs from the point
$p(0)=\left(34/3,109/93,-296730/604469+1/1813407\sqrt{1587984839746}\right)\simeq$
$ \left(11.33333\right.$, $1.17204$, $\left. 0.20401\right)\in{\cal G}$ to the
fixed point of $F^2,$ $p(1)\simeq\left(1.04794, 38.08255, 1.04794\right)$, over
the level surface
$${\textstyle V_1=\bar k:= {\frac {101}{272452149}}\,{\frac { \left( 923217+\sqrt
{1587984839746}
 \right)  \left( 69592724+3\,\sqrt {1587984839746} \right) }{-890190+
\sqrt {1587984839746}}}
 \simeq 0.24956.}$$

\begin{center}

\noindent {\small
\begin{tabular}{|c|c|c|c|c|}
\hline
 $t$ & Point $p$ & $T(p)$ & $\tau_{F^2}(p) $  & $\rho_{F^2}(p) $  \\
\hline \hline
$0  $ & $\left(11.33333, 1.17204, 0.20402\right) $ & $0.12622 $ & $0.03123 $ & $0.24738$ \\
\hline
$0.1  $ & $\left(10.30479, 1.19106, 0.18573\right) $ & $0.12620 $ & $0.03122 $ & $0.24739 $ \\
\hline
$0.2  $ & $\left(9.27625, 1.21480, 0.16809\right) $ & $0.12612 $ & $0.03120 $ & $0.24740  $ \\
\hline
$0.3  $ & $\left(8.24772, 1.24529, 0.15108\right) $ & $0.12596 $ & $0.03116 $ & $0.24742  $ \\
\hline
$0.4  $ & $\left(7.21918, 1.28585, 0.13473\right) $ & $0.12569 $ & $0.03110 $ & $0.24744  $ \\
\hline
$0.5  $ & $\left(6.19064, 1.34250, 0.11911\right) $ & $0.12526$ & $0.03010 $ & $0.24748 $ \\
\hline
$0.6  $ & $\left(5.16210, 1.42714, 0.10435\right) $ & $0.12455 $ & $0.03083 $ & $0.24755 $ \\
\hline
$0.7  $ & $\left(4.13356, 1.56733, 0.09081\right) $ & $0.12335 $ & $0.03055 $ & $0.24765 $ \\
\hline
$0.8  $ & $\left(3.10502, 1.84454, 0.07963\right) $ & $0.12107 $ & $0.03001 $ & $0.24784  $ \\
\hline
$0.9  $ & $\left(2.07648, 2.65147, 0.07637\right) $ & $0.11551 $ & $0.02867 $ & $0.24824  $ \\
\hline
$0.95  $ & $\left(1.56221, 4.16212, 0.09117\right) $ & $0.10807 $ & $0.02687 $ & $0.24866$ \\
\hline
$0.99  $ & $\left(1.15080, 12.78937, 0.23101\right) $ & $0.08942 $ & $0.02229 $ & $0.24932  $ \\
\hline
$0.999  $ & $\left(1.05822, 31.53212, 0.74344\right) $ & $0.07887 $ & $0.01968 $ & $0.24955 $ \\
\hline
$0.9999  $ & $\left(1.04897, 37.30369, 1.00598\right) $ & $0.07829 $ & $0.01954 $ & $0.24956  $ \\
\hline
$1 $ & $\left(1.04794, 38.08255, 1.04794\right) $ & $-- $ & $-- $ & $0.24957 $ \\
\hline
\end{tabular}}
\end{center}
\centerline{Table 4. Rotation number on $\{V_1=\bar k\simeq 0.24956\}$ when
$a=7/9.$}


\newpage

\appendix
\section*{Appendices}

\section{Proof of Lemma \ref{corbes}}\label{propoa}

 To  describe the foliation of $Q^+$,
induced by $\Delta(x,y;a,k)=0$ obtained for a fixed value of $a>0$, and varying
$k\geq k_c$, we solve the quadratic equation (with respect $k$):
$\Delta(x,y;a,k)=x^2y^2\,k^2+p_1(x,y;a)\, k+p_0(x,y;a)=0$. Thus the curve
$\Delta(x,y;a,k)=0$ in $Q^+$ can  also be described by two equations
$$k=m_\pm(x,y;a)=\dps{{\frac { \left( x+y+a+1\pm 2\,\sqrt {x+y+a} \right)
\left( x+1 \right)
 \left( y+1 \right) }{xy}}
 }.$$

We make the following  claims:

\vspace{0.5cm}

\rec{\bf Claim 1:} {\it For  any  fixed $k\geq k_c$ the following
statements hold

 (i) If $a>1$, there exist two values
$x_{1,k}<x_{2,k}$ such that equation (with unknown $y$)
\begin{equation}\label{kmenys}
  k=m_-(\bx,y;a),
\end{equation}
has two solutions if $\bx\in(x_{1,k},x_{2,k})$, one solution if
$\bx=x_{i,k}$ $i=1,2$, and none solution if $\bx\notin
[x_{1,k},x_{2,k}]$. This means that varying $x>0$, equation
(\ref{kmenys}) describes an oval $\zeta_k$.

(ii) If $a<1$, then there exist a value $x_{k}>1-a$ such that
equation (\ref{kmenys}) has two solutions if $\bx<x_{k}$, one
solution if $\bx=x_{k}$, and none solution if $\bx>x_k$. This
means that varying $x>0$, equation (\ref{kmenys}) describes a
curve consisting of a point $(x_k,y_k)$ and (from right to left)
two positive branches $y_1(x)<y_2(x)$ defined only for
$x\in(0,x_k)$. A more accurate analysis  will show that these two
branches meet at the point $(x,y)=(0,1-a)$. Therefore they
describe an oval, namely $\zeta_k$.}

\vspace{0.5cm}

\rec{\bf Claim 2:} {\it (i) For $k>k_c$, there exist two values
$x_{1,k}<x_c<x_{2,k}$ such that equation
\begin{equation}\label{kmes}
  k=m_+(\bx,y;a),
\end{equation} has two
solutions if $\bx\in(x_{1,k},x_{2,k})$, one solution if
$\bx=x_{i,k}$ $i=1,2$, and none solution if $\bx\notin
[x_{1,k},x_{2,k}]$. This means that varying $x>0$, equation
(\ref{kmes}) describes an oval $\gamma_k$.

(ii) The equation $k_c=m_{+}(\bx,y;a)$ has a unique solution if
$\bx=x_c$ and none solution if $\bx\neq x_c$.}

\vspace{0.5cm}

Since $m_+(x,y;a)>m_-(x,y;a)$  it is easy to see that each oval
$\gamma_k$ surrounds the corresponding oval $\zeta_k.$ From this
fact and the above claims the proof of the lemma follows.

Before proving both claims we establish some common facts. We fix $\bx>0$ and
we use the following notation:
$$
\dps{\frac{\partial m\pm}{\partial
y}}(\bx,y;a)=\dps{-(x+1)\left[\frac{\pm
f(\bx,y;a)+g(\bx,y;a)\sqrt{h(\bx,y;a)}}{\bx
y^2\sqrt{h(\bx,y;a)}}\right]},
$$
where $
  f(\bx ,y;a)= -y^2+y+2\bx+2a,$ $   g(\bx ,y;a)= -y^2+\bx+a+1$ and $
  h(\bx ,y;a)=\bx+y+a.$ The solutions in $Q^+$ of $\dps{\frac{\partial m\pm}{\partial y}}(\bx,y;a)=0$
are described by
\begin{equation}\label{fgh}\begin{array}{l}
\left(f^2-g^2\,h\right)(\bx,y;a)=(y+\bx+a-1)(-y^2-y+\bx+a-1)(-y^2+y+\bx+a).\end{array}\end{equation}
So this equation gives the local extrema of $y\to m_{\pm}(\bx,y;a)$.

It can be easily proved that
\begin{equation}\label{aaalinfinit}\lim_{y\to 0^+} m_{\pm}(\bx,y;a)=+\infty\mbox{ and } \lim_{y\to
+\infty} m_{\pm}(\bx,y;a)=+\infty,\end{equation} for all $\bx>0$.

Let us now proceed with the proof of both claims.

\vspace{0.5cm}

\rec{\it Proof of Claim 1.} (i) If $a>1$, and taking into account that $\bx>0$,
it is easy to see that $\dps{\frac{\partial m_{-}}{\partial y}}(\bx,y;a)=0$ if
and only if
\begin{equation}\label{simetria}
Q(\bx,y)=-y^2-y+\bx+a-1=0
\end{equation}
which has the unique positive solution $y=\ymin(\bx)=(-1+\sqrt{-3+4a+4\bx})/2.$

We point out that this solution is well defined for all $\bx>0$
and $a>1$, and that $y_{\min}(\bar{x})>0$. Hence, taking into
account equation (\ref{aaalinfinit}), the function $y\to
m_{-}(\bx,y;a)$ takes a minimum at the point $\ymin(\bx).$ So we
have that for each $\bx>0$ the functions $y\to m_{-}(\bx,y;a)$ are
decreasing in the interval $y\in(0,\ymin)$ and increasing in
$y\in(\ymin,+\infty)$.

Now we study the function $x\to m_{-}(x,\ymin(x);a)$. We have the following
facts:

\rec (I) $\lim\limits_{x\to 0^+} m_{-}(x,y_{\min}(x);a) =+\infty,$

\rec (II)  It is easy to see that at infinity
$m_{-}(x,y_{\min}(x);a)\sim x$, thus $$\lim\limits_{x\to +\infty}
m_{-}(x,y_{\min}(x);a)=+\infty.$$ \rec (III) Since the only
positive solutions of $\dps{\frac{\partial m_{-}}{\partial
x}}(x,y;a)=0$ are given by the solutions of the equation
$Q(y,x)=-x^2-x+y+a-1=0$, where $Q(x,y)$ is defined in
(\ref{simetria}), and $Q(x,y)=Q(y,x)=0$ if and only if
$x=y=x_*:=\sqrt{a-1},$ we obtain that $\dps{\frac{\partial
m_{-}}{\partial x}}(x,y_{\min}(x);a)\neq 0$  for $x\neq x_*,$ and
$\dps{\frac{\partial m_{-}}{\partial x}}(x_*,x_*;a)=0.$ This means
that $x_*$ is the unique critical point of the function $x\to
m_{-}(x,\ymin(x);a)$ when  $x>0$.

\rec (IV) It can be easily checked that $k_*=m_{-}(x_*,x_*;a)=(1+\sqrt{a-1})^2$
and  $ k_*< k_c.$

 Collecting the results in  items (I)--(IV) we obtain that
 the function
 $x\to m_{-}(x,y_{\min}(x);a),$
which gives the minimum values of each function $y\to
m_{-}(x,y;a),$ has a unique minimum at $x=x_*$, decreases from
$+\infty$ to $k_*$ for $x\in(0,x_*)$, and increases from $k_*$ to
$+\infty$ for $x\in (x_*,+\infty)$.

This proves that for any fixed $k\geq k_c$ there always exist two solutions
$x_{i,k}$ , $i=1,2$ of equation $ k=m_{-}(x,y_{\min}(x);a),$ see also Figure 3.

 For these two values, the
minimum of the functions $y\to m_{-}(x_{i,k},y;a)$ is $k$. This
implies that equation (\ref{kmenys}) only has one solution for
$\bx=x_{i,k}$ , $i=1,2$. Now observe that for all
$\bx\in(x_{1,k},x_{2,k})$, since the minimum values of the
functions $y\to m_{-}(\bx,y;a)$  are below $k$, we can conclude
that for these values of $\bx$ equation (\ref{kmenys}) only has
two solutions. Finally for $\bx\notin(x_{1,k},x_{2,k})$, since the
minimum values of the functions $y\to m(\bx,y;a)$ are greater than
$k$, equation (\ref{kmenys}) has no solutions. In summary,
equation (\ref{kmenys}) describes one and only one oval $\zeta_k$.



\centerline{\includegraphics[scale=0.65]{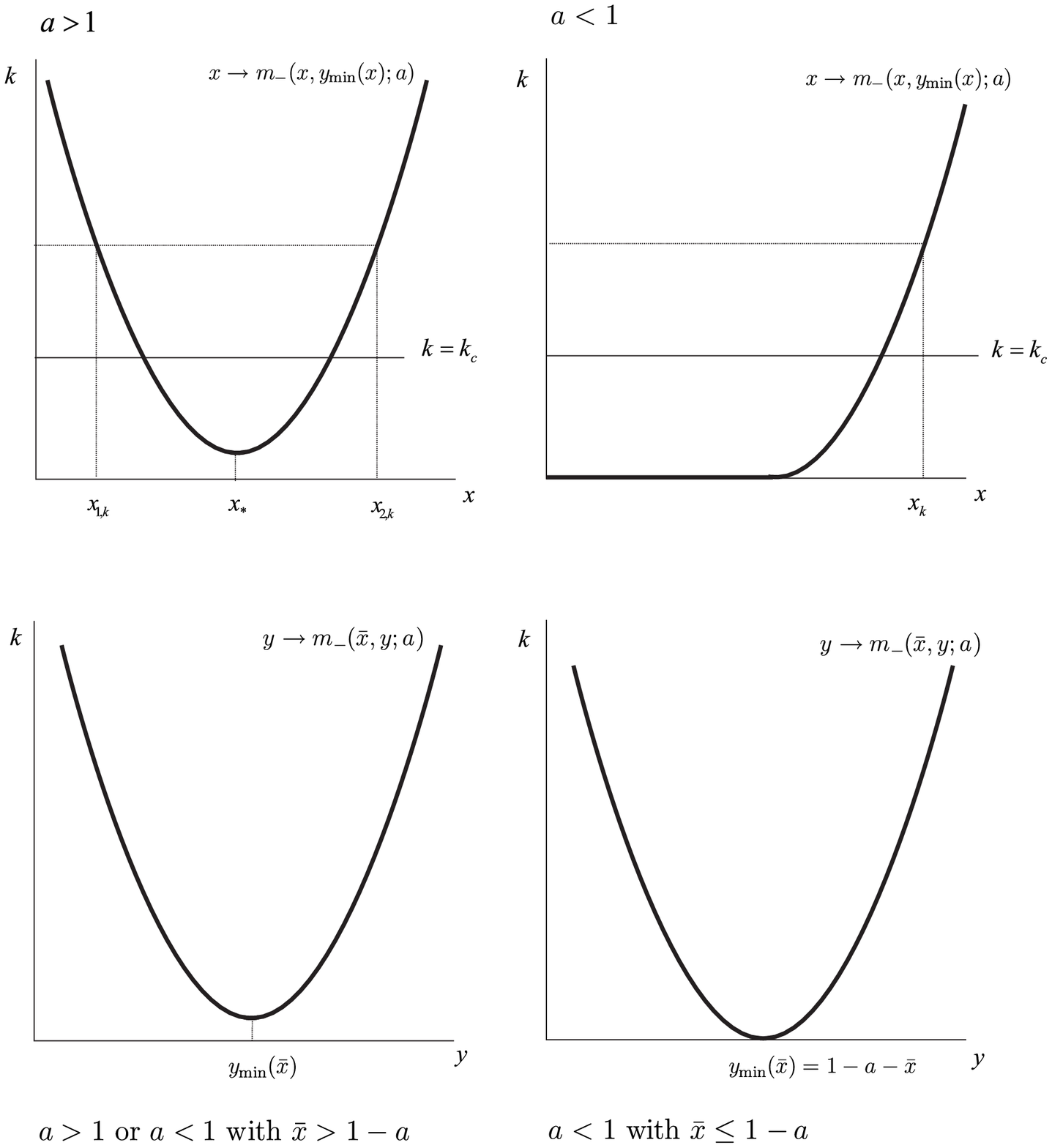}}
\centerline{Figure 3. Plot of some functions involved in the proof
of Lemma \ref{corbes}.}

\vspace{0.5cm}

(ii) If $a<1$, $\dps{\frac{\partial m_{-}}{\partial
y}}(\bx,y;a)=0$ at the curve  $y=\ymin(\bx)$, described by the
unique positive root of $(y+\bx+a-1)(-y^2-y+\bx+a-1).$ Thus
$$
\ymin(\bx)=\cases{
\begin{array}{ll}
  1-a-\bx & \mbox{if } \bx\leq 1-a, \\
  (-1+\sqrt{-3+4a+4\bx})/2 & \mbox{if } \bx>1-a,
\end{array}
}
$$
which is well defined for all $\bx>0$. Taking into account
equation (\ref{aaalinfinit}), $\ymin(x)$ gives a minimum for $y\to
m_{-}(\bx,y;a)$.

Now, we look at the function $x\to m_{-}(x,\ymin(x);a)$. First notice that for
$x\leq 1-a$, $m_{-}(x,\ymin(x);a)=m_{-}(x,x+a-1;a)=0$. In the region $x>1-a$,
we have
$$
m_{-}(x,\ymin(x);a)=\frac{n(x)(x+1)(1+\sqrt{-3+4a+4x})^2}{8x(x-1+a)},
$$
where $n(x)=\tilde{l}(s(x))$ with $s(x)=2(x+a)+\sqrt{-3+4x+4a}$, and
$\tilde{l}(s)=s+1+2\sqrt{2}\sqrt{s-1}$.

Observe that $x>1-a$ if and only if $s>3$, and that the function
$\tilde{l}(s)$ is monotonic increasing  from $0$ to $+\infty$ for
$s>3$. Hence for $x>1-a$, we have $m_{-}(x,\ymin(x);a)$ is
monotonic increasing from $0$ to $+\infty$.

Therefore, for $\bx<x_k$, the minimum of the functions  $y\to m_{-}(\bx,y;a)$
is  always below $k$, hence we can conclude that for these values of $\bx$
equation (\ref{kmenys}) has exactly two solutions, giving rise to two positive
branches $y_1(\bx)$ and $y_2(\bx)$. Since $\Delta(0,y;a,k)=(y+1)^2(a-1+y)^2$ we
can conclude that these two branches meet at the point $(0,1-a)$, which is an
order two contact point of $\Gamma_k$ with $\{x=0\}$. The symmetry of $\Delta$
with respect the line $\{y=x\}$, gives the other contact point $(1-a,0).$

For $\bx=x_k$ the minimum of the functions $y\to m_{-}(\bx,y;a)$ is $k$. This
implies that equation (\ref{kmenys}) has only one solution, giving the point
where the previous mentioned two branches meet.

Finally, observe that for all $\bx>x_k$, since the minimum values of the
functions $y\to m_{-}(\bx,y;a)$  are over $k$, equation (\ref{kmenys}) has no
solutions.

In summary equation (\ref{kmenys}) describes one and only one oval $\zeta_k$,
which has two contact points with the boundary of $Q^+$ at $(0,1-a)$ and
$(1-a,0)$. This ends the proof of Claim~1.

\vspace{0.5cm}

\rec{\it Proof of Claim 2.} We have to study the function $m_+$.
It is easy to see that $y\to m_{+}(\bx,y;a)$ has a unique minimum
at $y=y_{\min}(\bx)=\frac{1+\sqrt{1+4a+4\bx}}{2},$ the unique
positive solution of the equation $C(x,y):=-y^2+y+\bx+a=0$, and
that $m_{+}(\bx,y;a)$ is decreasing from infinity to
$m_{+}(\bx,\ymin(\bx);a)$ for $y<y_{\min}(\bx)$ and increasing to
infinity for $y>y_{\min}(\bx)$.

We need now to study the function  $x\to m_{+}(x,y_{\min}(x);a)$.
Observe that:

\rec (I)  It holds that
  $y_{\min}(x_c)=x_c$.

\rec (II) Since  $\lim\limits_{x\to 0^+} y_{\min}(x)=(1+\sqrt{1+4a})/2> 0$, we
get
   $\lim\limits_{x\to 0^+} m_{+}(x,y_{\min}(x);a) =+\infty.$

\rec (III) It is easy to check that at infinity
$m_{+}(x,y_{\min}(x);a)\sim x$, thus $$\lim\limits_{x\to +\infty}
m_{+}(x,y_{\min}(x);a)=+\infty.$$ \rec (IV) Since the only
positive solutions of $\dps{\frac{\partial m_{+}}{\partial
x}}(x,y;a)=0$ are given by the
 equation $C(y,x)=-x^2+x+y+a=0$, and
$ C(x,y)=  C(y,x)=0$ if and only if $ x=y=x_c=1+\sqrt{1+a},$ we
have that $\dps{\frac{\partial m_{+}}{\partial
x}}(x,y_{\min}(x);a)\neq 0$ for $x\neq x_c, $  and
$\dps{\frac{\partial m_{+}}{\partial x}}(x_c,y_{\min}(x_c);a)=0.$

Collecting the information summarized in (I)--(IV), we obtain that
 the function $x\to m_{+}(x,y_{\min}(x);a),$ which
gives the minimum values of each function $y\to m_{+}(x,y;a),$ has a unique
minimum at $x=x_c$, and decreases from $+\infty$ to $k_c$ for $x\in(0,x_c)$,
and increases from $k_c$ to $+\infty$ for $x\in (x_c,+\infty)$. A simple
computation, omitted here, shows that
$m_{+}(x_c,y_{\min}(x_c);a)=m_{+}(x_c,x_c;a)=k_c.$

This proves that for any fixed $k>k_c$ there exists only two solutions
$x_{i,k},$ $i=1,2$ of equation
\begin{equation}\label{lacorbadelsminims}
  k=m_{+}(x,y_{\min}(x);a),
  \end{equation}
such that $x_c\in(x_{1,k},x_{2,k}),$  see again Figure 3. For these two values,
the minimum of the functions $y\to m_+(x_{i,k},y;a)$ is $k$.  This means that
equation (\ref{kmes}) has only one solution for $\bx=x_{i,k}$ $i=1,2$. We note
that for all $\bx\in(x_{1,k},x_{2,k})$, since the minimum values of the
functions $y\to m_+(\bx,y;a)$  are below $k$, we can conclude that for these
values of $\bx$ equation (\ref{kmes}) has only two solutions. Finally for
$\bx\notin(x_{1,k},x_{2,k})$, since the minimum values of the functions $y\to
m_+(\bx,y;a)$ are greater than $k$,  equation (\ref{kmes}) has no solutions. In
summary equation (\ref{kmes}) describes one and only one oval $\gamma_k$. This
ends the proof of (i).

The above analysis of $x\to m_{+}(x,y_{\min}(x);a)$ shows that
$k_c=m(\bx,y;a)$ if and only if $\bx=y=x_c$. So the branch of
$\Gamma_{k_c}$ described by equation (\ref{kmes}) collapses to the
point $(x_c,x_c)$. This ends the proof of (i), and so the proof of
Claim 2.\qed

\section{Proof of Lemma \ref{corbes2}}\label{propob}

To describe the foliation of $Q^+$, induced by $\Delta(x,y;a,h)=0$ obtained for
a fixed value of $a>0$, and varying $h>h_c$, we can rewrite
$\Delta(x,y;a,h)=x^2y^2\,h^2+p_1(x,y;a)\, h+p_0(x,y;a)$ where
$$\begin{array}{ll}
  p_0(x,y;a)=&{x}^{2}{y}^{4}+ \left(2{x}^{3}+2{x}^{2}+ \left( -2a+2 \right) x  \right) {y}^{3
}\\
&+ \left({x}^{4}+2{x}^{3}+ \left( -4a+5 \right) {x}^{2}+ \left(
-4a+4 \right) x+{a}^{2}-2a+1
 \right) {y}^{2} \\
 &+\left(\left( -2a+2 \right) {x}^{3}+ \left( -4a+4 \right) {x}^{2}+ \left( 4+2{a}^{2}-6a \right) x+2+2{a}^{2}-4a
 \right) y\\
 &+\left( {a}^{2}-2a+1 \right) {x}^{2}
 + \left( 2+2{a}^{2}-4a \right) x+{a}^{2}-2a+1, \mbox{ and}\\
  p_1(x,y;a)=&((-2x^2-4x))y^3+(-2x^3-10x^2+(-2a-6)x)y^2\\&+(-4{x}^{3}+ \left( -2a-6 \right) {x}^{2}+
  \left( -2-2a \right) x)y.
\end{array}$$
Thus the curve $\Delta(x,y;a,h)=0$ in $Q^+$ can  also be described
by two functions $$h=m_\pm(x,y;a)=\dps{{\frac { \left(
yx+2x+2y+a+1\pm 2\sqrt {d(x,y;a)} \right)  \left( 1+x+y \right)
}{yx}} },$$ where$$d(x,y;a)=
{x}^{2}y+x{y}^{2}+{x}^{2}+{y}^{2}+(a+2)xy +(a+1)x+(a+1)y+ a.$$

 As in the previous appendix, to prove the lemma we make two claims:

\vspace{0.5cm}

\rec{\bf Claim 1:} {\it For $h\geq h_c$ and for all fixed $\bx>0$, we will see
that there exist two solutions $y_1(\bx)<y_2(\bx)$ of the equation
\begin{equation}\label{kmenys2}
  h=m_{-}(\bx,y;a),
\end{equation}
that give rise to the two branches of $\Gamma_h$, $y_1(x)$ and $y_2(x)$ given
in the statement of the lemma. Moreover, $\lim\limits_{x\to 0^+}
y_i(x)=+\infty$ and $\lim\limits_{x\to +\infty} y_i(x)=0^+.$}

\vspace{0.5cm}

\rec{\bf Claim 2:} {\it (i) For $h>h_c,$ if we consider the equation
\begin{equation}\label{kmes2}
  h=m_+(\bx,y;a),
\end{equation}
then there exist two values $x_{1,h}<x_{2,h}$ such that equation (\ref{kmes2})
has two solutions if $\bx\in(x_{1,h},x_{2,h})$, one solution if $\bx=x_{i,h}$
$i=1,2$, and none solution if $\bx\notin [x_{1,h},x_{2,h}]$. This means that
varying $x>0$, equation (\ref{kmes2}) describes the oval $\gamma_h$.

(ii) The equation $h_c=m_{+}(\bx,y;a)$ has a unique solution if $\bx=x_c$ and
none solution if $\bx\neq x_c$.}

Note that, since $m_+(x,y;a)>m_-(x,y;a),$ the solutions of the equation
(\ref{kmes2}) (whenever they exist) are contained in the interval
$(y_1(x),y_2(x))$ defined by equation (\ref{kmenys2}). This implies that the
oval $\gamma_h$ is contained between the two branches $y_1(x)$ and $y_2(x)$.
So, by using the above two claims,  the lemma follows.

Before giving the proof of the above claims we
 establish some common facts. We fix $\bx>0$ and we use the following
notation:
$$
\dps{\frac{\partial m\pm}{\partial y}}(\bx,y;a)=\dps{\frac{\mp
f(\bx,y;a)+g(\bx,y;a)\sqrt{h(\bx,y;a)}}{\bx
y^2\sqrt{h(\bx,y;a)}}},
$$
where $$
\begin{array}{rl}
  f(\bx ,y;a)=& \left(-2\bx  -2\right) {y}^{3}+ \left(-{\bx }^{2}+ \left( -2-a \right) \bx -1-a  \right) {y}^{2}\\
  &+ \left({\bx }^{3}+ \left( 3+a \right) {\bx }^{2}+ \left( 2a+3 \right) \bx +1+a  \right) y\\
  &+2{\bx }^{3}+ \left( 2a+4 \right) {\bx }^{2}+ \left( 2+4a \right) \bx +2
a, \\
   g(\bx ,y;a)=& \left( \bx +2 \right) {y}^{2}-2{\bx }^{2}+ \left( -a-3 \right) \bx -1-a, \\
  h(\bx ,y;a)=&\left( y+1 \right)  \left( \bx+1 \right)  \left( \bx+a+y \right).
\end{array}
$$
 So $\dps{\frac{\partial m\pm}{\partial y}}(\bx,y;a)=0$ if and only if
$$\begin{array}{l} f^2(\bx,y;a)-g^2(\bx,y;a)\,{h(\bx,y;a)}=(\bx+1)\left(-y^3-y^2+(1+\bx)y+(1+\bx)(a+\bx)\right)\cdot\\
\cdot(\bx y-a+1)\left(\bx
y^2+(1+\bx)(a-1+\bx)y+(1+\bx)(a-1)\right)=0.\end{array}$$ This
equation gives the local extrema of $y\to m_{-}(\bx,y;a)$ and
$y\to m_{+}(\bx,y;a)$.

It is not difficult to see  that
\begin{equation}\label{aaalinfinit2}\lim_{y\to 0^+} m_{\pm}(\bx,y;a)=+\infty\quad \mbox{ and } \lim_{y\to
+\infty} m_{\pm}(\bx,y;a)=+\infty,\end{equation} for all $\bx>0$.

\vspace{0.5cm}

\rec{\it Proof of Claim 1.} If $a>1$, taking into account that
$\bx>0$, it is easy to see that $\dps{\frac{\partial
m_{-}}{\partial y}}(\bx,y;a)=0$ if and only if $y_0=(a-1)/\bx.$
Furthermore $m_{-}(\bx,(a-1)/\bx;a)= 0$ for all $x>0$, which is a
minimum of $y\to m_{-}(\bx,y;a)$. Taking into account equation
(\ref{aaalinfinit2}), we have that for each $\bx>0$ the functions
$y\to m_{-}(\bx,y;a)$ are decreasing from $+\infty$ to $0$ for
$y\in(0,y_0)$ and increasing from $0$ to $+\infty$ if
$y\in(y_0,+\infty)$. This proves that in this case equation
(\ref{kmenys2}) always has two solutions for any $h>0$, in
particular for any $h\geq h_c$.

If $0<a<1$, $\dps{\frac{\partial m_{-}}{\partial y}}(\bx,y;a)=0$ if and only if
$y=y_q(\bx)$, where $y_q(\bx)$ is the only positive solution of the quadratic
equation $\bx y^2+(1+\bx)(a-1+\bx)y+(1+\bx)(a-1)=0$ which taking into account
(\ref{aaalinfinit2}) gives  a minimum of $y\to m_{-}(\bx,y;a)$. To see that
equation (\ref{kmenys2}) always has two solutions for $h\geq h_c$, we only have
to see that $m_{-}(x,y_q(x);a)<h_c$ for all $x>0$. Since $ m_{-}(x,y;a)=1-a$ if
and only if $x y^2+(1+x)(a-1+x)y+(1+x)(a-1)=0,$ we have that
$m_{-}(x,y_q(x);a)= 1-a.$ On the other hand it  is not difficult to check that
$1-a<h_c.$ So the first part of  the claim is proved.

To end the proof it remains to see that $\lim\limits_{x\to 0^+} y_i(x)=+\infty$
and $\lim\limits_{x\to +\infty} y_i(x)=0^+,$ where $y_1(x)$ and  $y_2(x)$ are
the two branches of $\Delta=0.$ We observe that each curve $\Delta(x,y;a,h)=0$
is symmetric with respect the axis $y=x$. So it is equivalent to see that
$\lim\limits_{x\to+\infty} y_i(x)=0$ , $i=1,2$ or $\lim\limits_{x\to 0^+}
y_i(x)=+\infty$ $i=1,2$. So we will prove the first equality. To do this we
will study if there arrive any branch of $\Gamma_h$ to the ``infinity line'' in
the projective space $\pr$ in the direction ${y=0}$.

For each affine curve  $\Delta(x,y;a,h)=0$ we can consider the
projectivized curve in $\pr$ (given in homogeneous coordinates
$(x,y,u)$) $ \tilde{\Gamma}_h=\{\tD(x,y,u;a,h)=0\}$, where

\rec $ \tD(x,y,u;a,h)=\left( {a}^{2}-2a+1 \right) {u}^{6}+ \left(
2+2{a}^{2}-4a \right) x{u}^{5}+ \left( {a}^{2}-2a+1 \right)
{x}^{2}{u}^{4}+\\
 \left( 2+2{a}^{2}-4a \right) y{u}^{5}+ \left( 4+2{a}^{2}-6a-2
h-2ka \right) xy{u}^{4}+ \left( -4a-2ka+4-6h \right) {x}^{2} y{u}^{3}+\\ \left(
2-4h-2a \right) {x}^{3}y{u}^{2}+ \left( {a}^{2}-2 a+1 \right) {y}^{2}{u}^{4}+
\left( -4a-2ka+4-6h \right) x{y}^{ 2}{u}^{3}\\+ \left( {h}^{2}-4a-10h+5 \right)
{x}^{2}{y}^{2}{u}^{2}+
 \left( -2h+2 \right) {x}^{3}{y}^{2}u+{y}^{2}{x}^{4}+ \left( 2-4h-
2a \right) x{y}^{3}{u}^{2}\\+ \left( -2h+2 \right)
{x}^{2}{y}^{3}u+2 {y}^{3}{x}^{3}+{y}^{4}{x}^{2}$. In the local
chart $\{x\neq 0\}$ this curve is given by (just taking $x=1$)
$\tilde{\Gamma}_h=\{\tD(1,y,u;a,h)=0\}$, where:

\rec $\tD(1,y,u;a,h)=\left( {a}^{2}-2a+1 \right) {u}^{6}+
\left( 2+2{a}^{2}-4a \right) {u}^{5}+ \left( {a}^{2}-2a+1 \right) {u}^{4}+\\
\left( 2+2{ a}^{2}-4a \right) y{u}^{5}+ \left( 4+2{a}^{2}-6a-2h-2ka
 \right) y{u}^{4}+ \left( -4a-2ka+4-6h \right) y{u}^{3}+ \\\left(
2-4h-2a \right) y{u}^{2}+ \left( {a}^{2}-2a+1 \right) {y}^{2}{u} ^{4}+ \left(
-4a-2ka+4-6h \right) {y}^{2}{u}^{3}+\\ \left( {h}^{2} -4a-10h+5 \right)
{y}^{2}{u}^{2}+ \left( -2h+2 \right) {y}^{2}u+ {y}^{2}+ \left( 2-4h-2a \right)
{y}^{3}{u}^{2}+ \left( -2h+2
 \right) {y}^{3}u+2{y}^{3}+{y}^{4}.
$

We want to prove is that there are two ``affine'' branches of
$\tilde{\Gamma}_h$ arriving at the point of the infinity line with
coordinates  $(y,u)=(0,0)\in \tilde{\Gamma}_h$. Since
$\tD(1,y,u;a,h)=y^2+2y^3+(2-2h)y^2u+2(1-a-2h)y u^2+O((x,y)^4)$, we
need to perform the blow--up: $(y,u)=(v\,u,u)$, which after
removing the factor $u^2$ transforms $\tilde{\Gamma}_h$ into

\rec $\tilde{\Gamma}_h^*=\{
 \left[ \left( {a}^{2}-2\,a+1 \right) {v}^{2}+ \left( 2+2\,{a}^{2}-4\,a \right) v+{a}^{2}-2\,a+1
 \right] {u}^{4}\\+ \left[ \left( 2-4\,h-2\,a \right) {v}^{3}\\+ \left( -4\,a-2\,ka+4-6\,h
 \right) {v}^{2}+ \left( 4+2\,{a}^{2}-6\,a-2\,h-2\,ka \right) v+2+\right.\\\left. 2\,{
a}^{2}-4\,a
 \right] {u}^{3}+ \left[{v}^{4}+ \left( -2\,h+2 \right) {v}^{3}+ \left( {h}^{2}-4\,a-10\,h+5
 \right) {v}^{2}+ \left( -4\,a-2\,ka+4-6\,h \right) v+\right.\\\left.{a}^{2}-2\,a+1
 \right] {u}^{2}+ \left[2\,{v}^{3}+ \left( -2\,h+2 \right) {v}^{2}+ \left( 2-4\,h-2\,a
 \right) v \right]u+{v}^{2}
 =0\}.$

The intersection $\tilde{\Gamma}_h^*$ with $\{v=0\}$  are the
points $(v,u)=(0,0)$, and $(v,u)=(0,-1)$. But the last point is
not interesting for us since, the affine region $\{(x,y): x>0,
y\geq 0\}$ in this local coordinates corresponds with $\{(v,u):
v\geq 0, u>0\}$. The directions of approach of
$\tilde{\Gamma}_h^*$ to $(v,u)=(0,0)$ are given by:
$u=\lambda_{\pm}\,v$, where
$$\lambda_{\pm}=\frac{(2h+a-1)\pm 2\sqrt{h(h+a-1)}}{(a-1)^2}.$$
It is easy to check that for $a\neq 1$, both $\lambda_{\pm}$ are
positive, therefore there exist two branches of
$\tilde{\Gamma}_h^*$ arriving at the singular point $(0,0)$ in
$\{(v,u): v\geq 0, u>0\}$. As only the two branches of $\Gamma_h$
described by $y_i(x)$, $i=1,2$ are defined when $x\to+\infty$
these ones are the two branches described by the blow--up
procedure. This ends the proof of the claim.

\vspace{0.5cm}

\rec{\it Proof of Claim 2}. It is easy to see that $y\to m_{+}(\bx,y;a)$ has a
unique minimum at $y=y_c(\bx)$, where, by Descartes' rule,  $y_c(\bx)$ is the
only positive solution  of the equation
$C(x,y):=-y^3-y^2+(1+\bx)y+(1+\bx)(a+\bx)=0$, and its decreasing at
$y<y_c(\bx)$ at increasing for $y>y_c(\bx)$.

Now we state three facts concerning the curve $y=y_c(x)$, which
are relevant to the study of $x\to m_{+}(x,y_c(x);a)$.
\begin{itemize}
  \item[(a)] An straightforward computation shows that
  $y_c(x_c)=x_c$.
  \item[(b)] Once again, applying Descartes's Rule on the cubic $C(x,y)=0,$
  we have that $\lim\limits_{x\to 0^+} y_c(x)> 0$,
  which implies that  $\lim\limits_{x\to 0^+} m_{+}(x,y_c(x);a)=+\infty.$
  \item[(c)] A detailed analysis of the asymptotic expansion of  $y_c(x)$ at
  infinity gives that $y_c(x)\sim (\sqrt[3]{216}/6)\sqrt[3]{x^2}$,
  hence   $\lim\limits_{x\to +\infty} y_c(x)=+\infty$,  and as a consequence
  $\lim\limits_{x\to +\infty} m_{+}(x,y_c(x);a)=+\infty.$
\end{itemize}
Since the only positive solutions of $\dps{\frac{\partial m_{+}}{\partial
x}}(x,y;a)$ are given by the cubic equation
$C(y,x)=-x^3-x^2+(1+y)x+(1+y)(a+y)=0$, and $ C(x,y)=
  C(y,x)=0$ if and only if $x=x_c$  and $y=x_c,$
we get that $\dps{\frac{\partial m_{+}}{\partial x}}(x,y_c(x);a)\neq 0$  for
$x\neq x_c,$  and $\dps{\frac{\partial m_{+}}{\partial x}}(x_x,y_c(x_c);a)=0.$

This means that the function $x\to m_{+}(x,y_c(x);a),$ which gives
the minimum values of each function $y\to m_{+}(x,y;a),$ has a
unique minimum at $x=x_c$,  decreases from $+\infty$ to $h_c$ for
$x\in(0,x_c)$, and increases from $h_c$ to $+\infty$ for $x\in
(x_c,+\infty)$. A simple computation shows that
$m_{+}(x_c,y_c(x_c);a)=m_{+}(x_c,x_c;a)=h_c.$

The above results  prove (see Figure 5), that for any fixed $h>h_c$ there
exists only two solutions  $x_{i,h},$  $i=1,2$ of equation
\begin{equation}\label{lacorbadelsminims2}
  h=m_{+}(x,y_c(x);a),
  \end{equation}
such that $x_c\in(x_{1,h},x_{2,h})$. For these two values,  the
minimum of the  functions $y\to m(x_{i,h},y;a)$ is $h$. This means
that equation (\ref{kmes2}) has only one solution for
$\bx=x_{i,h}$ $i=1,2$. Now we observe that for all
$\bx\in(x_{1,h},x_{2,h})$, since the minimum values of the
functions $y\to m(\bx,y;a)$  are below $h$, we can conclude that
for these values of $\bx$ equation (\ref{kmes2}) only has two
solutions. Finally for $\bx\notin(x_{1,h},x_{2,h})$, since the
minimum values of the functions $y\to m(\bx,y;a)$ are greater than
$h$, equation (\ref{kmes2}) has no solutions. In summary equation
(\ref{kmes2}) describes one and only one oval $\gamma_h$.

\vspace{0.5cm}


\centerline{\includegraphics[scale=0.8]{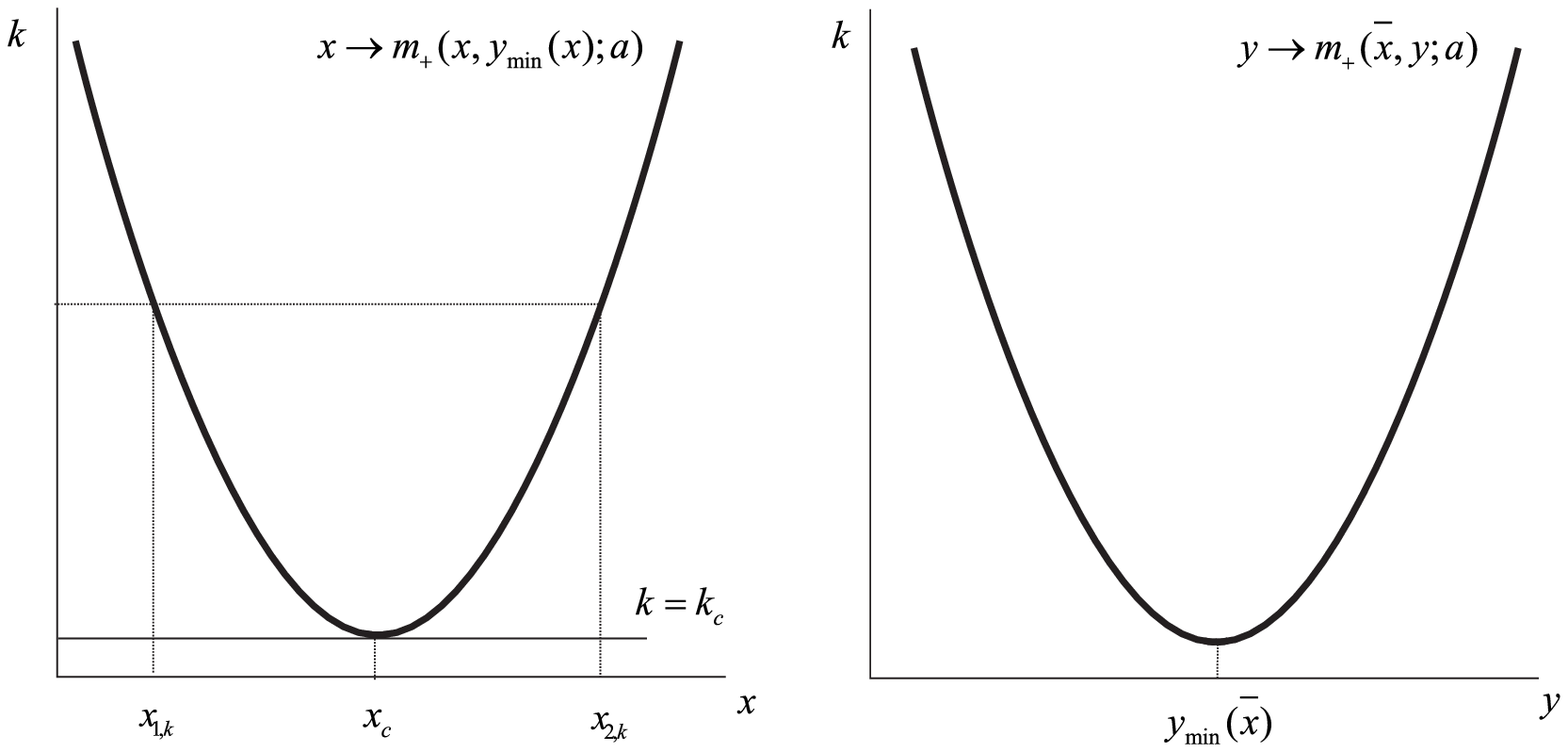}}
\centerline{Figure 5.  Plot of some functions involved in the
proof of Lemma \ref{corbes2}. }

\vspace{0.5cm}

The above analysis of $x\to m_{+}(x,y_c(x);a)$ shows that $h_c=m(\bx,y;a)$ if
and only if $\bx=y=x_c$. So the branch of $\Gamma_{h_c}$ described by equation
(\ref{kmes2}) is only the point $(x_c,x_c)$. This ends the proof of the claim.
\qed

\section{Proof of Proposition \ref{lem3}}\label{appd}

(i) Set
$$
h(x):=\left.{V_1}\right|_{\lin}=V_1\left(x,\dps{\frac{x+a}{x-1}},x\right)=
\dps{\frac{(2x+a-1)^2(x+1)^2}{x(x+a)(x-1)}},\mbox{
for } x>1.
$$
Trivial but tedious computations show that  the unique solution of
$h'(x)=0$ such that $x>1$ is $x=x_c$ which is a minimum, and that
$ \lim\limits_{x\rightarrow 1^+}h(x)=\lim\limits_{x\rightarrow
+\infty}h(x)=+\infty.$ Hence for any $k>k_c$, the equation
$V_1(x,(x+a)/(x-1),x)=k$ has a unique solution in $(1,x_c)$ and a
unique solution in $(x_c,+\infty)$. Hence $L_k\cap\lin$ consists
of two points, say $p_1$ and $p_2$. Recall that $\lin$ is the
curve of $2$--periodic points of $F$, hence $\{p_1,p_2\}$ is a
$2$--periodic orbit and (i) holds.

To prove statement (ii) we will see that if $k>k_c$ the locus of
non transversal intersections of $L_k$  with $\sg$  is the empty
set. Consider the system $G=0$ and
\begin{eqnarray}
  (V_1)_x/(V_1)_y & =&G_x/G_y,\label{eq1} \\
 (V_1)_x/(V_1)_z & =&G_x/G_z,\label{eq2}
\end{eqnarray}
obtained by imposing $\{\nabla V_1 \| \nabla G\}\cap\{G=0\}$.

The only positive solutions of equation (\ref{eq1}) are given by the zeroes of
$m_1:=xaz+xa{z}^{2}+2a{y}^{3}x+3{x}^{2}az+5{y}^{4}+3{y}^{5}+2{y}^{3}+3{x}^{2}a{z}^{2}+2{x}^{3}az+2{x}^{3}a{z}^{2}+x{y}^{2
}z-{y}^{2}x{z}^{2}+6a{y}^{2}z+4a{y}^{3}z+2x{y}^{3}z+2a{y}^{2}x
-3{y}^{2}{x}^{2}{z}^{2}+{x}^{2}z+x{z}^{2}+x{z}^{3}+3{x}^{3}z-2{y
}^{2}{x}^{3}{z}^{2}-4{y}^{2}{x}^{3}+2{x}^{3}{z}^{3}+2{x}^{4}{z}^
{2}+2{x}^{4}z-2{y}^{4}{x}^{2}-2{y}^{3}{x}^{3}+3{x}^{2}{z}^{3}+
4{x}^{2}{z}^{2}+5{x}^{3}{z}^{2}-4{y}^{2}{x}^{2}a-2{y}^{3}{x}^{
2}a-2{y}^{2}{x}^{3}z-2y{x}^{2}a+2ayz-2{y}^{3}{x}^{2}z-2{x}^{
2}yz-2y{x}^{3}+3{y}^{4}x+2{y}^{3}{z}^{2}+5a{y}^{4}+3{y}^{2}{
z}^{2}+3{y}^{2}z+3a{y}^{2}+8{y}^{3}z+4{y}^{3}x+8a{y}^{3}+3
{a}^{2}{y}^{2}-4{y}^{3}{x}^{2}-3{y}^{2}{x}^{2}-y{x}^{2}+{a}^{2}y+2
{a}^{2}{y}^{3}+{y}^{2}x+y{z}^{2}+5\,{y}^{4}z-7\,{y}^{2}{x}^{2}z,$

\rec and its zeroes over $\sg$ are given by the the zeroes of $r_1$, where it
satisfies $m_1=q_1\,G\,+r_1$, and it is given by $ r_1:=\left(
2xz-2{x}^{3}z+2x{z}^{2}-2{x}^{3}{z}^{2} \right) {y}^{2}+ \left(
4x{z}^{2}+2{x}^{2}z+2xa{z}^{2}-2{x}^{3}z\right.$ $\left.+4{x}^{2}{
z}^{2}+2x{z}^{3}+2{x}^{2}a{z}^{2}-2{x}^{3}{z}^{2}-2{x}^{4}{z}^
{2}-2{x}^{4}z+2{x}^{2}{z}^{3}+2xaz+2xz+2{x}^{2}az \right) y\,+$\break $
2{x}^{3}{z}^{3}+2xaz+2x{z}^{2}+2x{z}^{3}+4{x}^{2}a{z}^{2}+2
{x}^{3}a{z}^{2}+4{x}^{2}az+2{x}^{3}{z}^{2}+4{x}^{2}{z}^{3}+2
{x}^{3}az+2xa{z}^{2}+4{x}^{2}{z}^{2}.$

The only positive solutions of equation (\ref{eq2}) are given by the zeroes of
$m_2:=x-z$ and
$m_3:=2xaz+2xa{z}^{2}+2{x}^{2}az-{y}^{3}+2{x}^{2}a{z}^{2}-ayx+2xy{z}^{2}-2x{y}^{2}z-a{y}^{2}z-a{y}^{2}x+4{x}^{2}z+4x{z}^{
2}+3x{z}^{3}+3{x}^{3}z+2{x}^{3}{z}^{3}+5{x}^{2}{z}^{3}+9{x}^
{2}{z}^{2}+5{x}^{3}{z}^{2}-ya-ayz+2{x}^{2}yz+2y{z}^{2}{x}^{2}-xy
-yz-{y}^{2}-{y}^{2}{z}^{2}-2{y}^{2}z-a{y}^{2}-{y}^{3}z-{y}^{3}x-{y}^
{2}{x}^{2}-y{x}^{2}-2{y}^{2}x-y{z}^{2}+xz$.

Suppose that $m_2=0$ so $z=x$ and $r_1(x,y,x)=-2x^2(x+1)^2(1+x+y)(-a-x-y+xy)$,
hence the  positive solutions  of $r_1(x,y,x)=0$ are given by $y=(x+a)/(x-1)$,
thus the points in $\lin$. But as mentioned above $\lin\cap\sg=(x_c,x_c,x_c)$.

To study the zeroes of $m_3$ in $\sg$, we consider the zeroes of $r_2$,
satisfying  $m_2=q_2\,G\,+r_2$. But it is  given by $
r_2:=2xz(x+1)(z+1)(a+x+y+z+xz)$, and therefore $r_2=0$  has not positive
solutions.

In summary, if $k>k_c$, then  $\nabla V_1$ is never parallel to
$\nabla G$ over $\sg$ and hence $L_k\trans \sg$.

(iii) Recall that $\sg$ is defined by the equation
$$G=-y^3-(x+z+a+1)y^2-(x+z+a)y+xz(x+1)(z+1)=0.$$
By applying Descartes' Rule on $\sg$ we obtain that for all $x>0$
and $z>0$ there exist a unique $y(x,z)>0$ solution of $G=0$.
Consider  the function $ v(x,z):=V_1(x,y(x,z),z).$ Now the proof
is done in two steps: (I) The only singular point of $v$ is
$(x_c,x_c)$, which is a minimum.   (II) Each level curve
$v(x,z)=k>k_c$  is a closed curve surrounding $(x_c,x_c)$.

Step I: To find the singular points of $h(x,z)$ we look for the solutions of
system $$\cases{\begin{array}{l}
v_x=(V_1)_x+(V_1)_y\dps{\frac{\partial y}{\partial x}}=0,\\
v_z=(V_1)_z+(V_1)_y\dps{\frac{\partial y}{\partial z}}=0,
\end{array}}
$$
such that $x>0$ and $z>0$. The only factors in $v_x$ and $v_z$
giving rise to such solutions are

\rec $m:=xaz+xa{z}^{2}+2a{y}^{3}x+3{x}^{2}az+5{y}^{4}+3{y}^{5}+2
{y}^{3}+3{x}^{2}a{z}^{2}+2{x}^{3}az+2{x}^{3}a{z}^{2}+x{y}^{2
}z-{y}^{2}x{z}^{2}+6a{y}^{2}z+4a{y}^{3}z+2x{y}^{3}z+2a{y}^{2}x
-3{y}^{2}{x}^{2}{z}^{2}+{x}^{2}z+x{z}^{2}+x{z}^{3}-2{y}^{3}{x}^{2}
a-4{y}^{2}{x}^{2}a+3{x}^{3}z+2{x}^{3}{z}^{3}+2{x}^{4}{z}^{2}+2
{x}^{4}z+3{x}^{2}{z}^{3}+4{x}^{2}{z}^{2}+5{x}^{3}{z}^{2}-2{y
}^{2}{x}^{3}{z}^{2}-4{y}^{2}{x}^{3}+3{y}^{4}x+2{y}^{3}{z}^{2}+5
a{y}^{4}+3{y}^{2}{z}^{2}+3{y}^{2}z+3a{y}^{2}+8{y}^{3}z+4{y
}^{3}x+8a{y}^{3}+3{a}^{2}{y}^{2}-4{y}^{3}{x}^{2}-3{y}^{2}{x}^{
2}-2y{x}^{3}-y{x}^{2}+{a}^{2}y+2{a}^{2}{y}^{3}+{y}^{2}x+y{z}^{2}+5
{y}^{4}z-7{y}^{2}{x}^{2}z+2ayz-2y{x}^{2}a-2{y}^{2}{x}^{3}z-2
{y}^{3}{x}^{2}z-2{x}^{2}yz-2{y}^{4}{x}^{2}-2{y}^{3}{x}^{3}$

\rec and

\rec
$n:=xaz+3xa{z}^{2}+4a{y}^{3}x+{x}^{2}az+5{y}^{4}+3{y}^{5}+2{y}^{3}+3{x}^{2}a{z}^{2}+2xa{z}^{3}+2{x}^{2}a{z}^{3}+2{x}^
{2}{z}^{4}+2x{z}^{4}-2{y}^{2}{x}^{2}{z}^{3}-2{y}^{2}x{z}^{3}-2
y{z}^{2}a-4{y}^{2}{z}^{2}a-2{y}^{3}{z}^{2}a-2{y}^{3}x{z}^{2}+2
ayx-2xy{z}^{2}+x{y}^{2}z-7{y}^{2}x{z}^{2}+2a{y}^{2}z+2a{y}^{3}
z+2x{y}^{3}z+6a{y}^{2}x-3{y}^{2}{x}^{2}{z}^{2}+{x}^{2}z+x{z}^{2}
+3x{z}^{3}+{x}^{3}z+2{x}^{3}{z}^{3}+5{x}^{2}{z}^{3}+4{x}^{2}{z
}^{2}+3{x}^{3}{z}^{2}+5{y}^{4}x-4{y}^{3}{z}^{2}+5a{y}^{4}-3{
y}^{2}{z}^{2}+{y}^{2}z+3a{y}^{2}+4{y}^{3}z+8{y}^{3}x+8a{y}^{3}
+3{a}^{2}{y}^{2}+2{y}^{3}{x}^{2}+3{y}^{2}{x}^{2}+y{x}^{2}+{a}^{2
}y-4{y}^{2}{z}^{3}-2y{z}^{3}-2{y}^{4}{z}^{2}-2{y}^{3}{z}^{3}+2
{a}^{2}{y}^{3}+3{y}^{2}x-y{z}^{2}+3{y}^{4}z-{y}^{2}{x}^{2}z$,

\rec respectively. Here $y$ denotes $y(x,z).$

The common zeroes of $m$ and $n$ in $\sg$   are given by the zeroes of the
functions $r$ and $s$ respectively, where $m=pG+r$ and $n=qG+s$ for some
polynomials $p$ and $q.$ These functions are
$$
\begin{array}{l}
  r =2xz\,(x+1)\,(1+x+y)\,(z+1)\,(a+y+z-xy), \\
  s =-2xz\,(x+1)\,(1+y+z)\,(z+1)\,(-a-x-y+yz).
\end{array}
$$
The only positive solutions of $r=0$ and $s=0$ are given by
$(x,(x+a)/(x-1),x)$, which are the points of $\lin.$ Since $\lin\cap
\sg=(x_c,x_c,x_c),$ the proof of (I) is finished.

Step II: Note the following facts: from statement (ii), $\sg\trans
L_k;$ the level curves  $v(x,z)=k,$ are defined by $L_k\cap \sg$
and thus for analytic equations  and for $k>k_c$ they have no
critical points (since otherwise the hamiltonian vector field
$-v_z\partial_x+v_x\partial_z$ would have another critical point
than $(x_c,x_c)$  in contradiction with Step I); The sets $L_k$
are compact. From all them  we conclude that for a fixed $k>k_c$,
each level set  of $v(x,z)=k$  is diffeomorphic to a finite union
of closed curves. To prove that indeed  it is formed  by an unique
closed curve it suffices to show that the function $x\to v(x,x_c)$
is monotonic in $(x_c,+\infty).$ To see this  we will prove that
it has the unique critical point $x=x_c.$ To this end note that
the only positive solutions of $r=0$ are given by the factor
$a+y+z-xy$, hence $x=(a+y+z)/y$. Now
$G((a+y+x_c)/y,y,x_c)=(y-x_c)\,q_4(y)/y^2$, where $q_4$ is a
degree four polynomial in $y$ without positive solutions. Thus, as
we wanted to see, $x=x_c$ is the unique critical point of $x\to
v(x,x_c)$ and (iii) follows. \qed

\end{document}